\documentclass[11pt,twoside]{article}

\usepackage{graphicx}
\usepackage{multirow}
\usepackage{fancyhdr}
\usepackage{array}
\usepackage{amsfonts}
\usepackage{amsmath}
\usepackage{amssymb}
\usepackage{amsthm}

\usepackage{subfigure}

\usepackage{graphics} 
\usepackage{epsfig} 
\usepackage{latexsym}

\usepackage[textwidth=16cm,textheight=22cm,top=3cm,bottom=3cm]{geometry}
\usepackage{fancyhdr}

\newtheorem{lemma}{Lemma}

\newtheorem{remark}{Remark}
\theoremstyle{definition}

\usepackage{color}

\graphicspath{{./figures/}}

\newcommand{\drop}[1]{}
\newcommand{\no}{\noindent}
\newcommand{\fer}[1]{(\ref{#1})}
\newcommand{\qtext}[1]{\quad\text{#1}}
\newcommand{\qtextq}[1]{\quad\text{#1}\quad}

\newcommand{\peq}[1]{{\text{\tiny \emph{#1}}}}

\newcommand{\bi}{\mathbf{i}}
\newcommand{\bj}{\mathbf{j}}
\newcommand{\bJ}{\mathbf{J}}
\newcommand{\bk}{\mathbf{k}}

\newcommand{\bx}{\mathbf{x}}
\newcommand{\by}{\mathbf{y}}
\newcommand{\bz}{\mathbf{z}}

\newcommand{\bzero}{\mathbf{0}}

\newcommand{\cA}{\mathcal{A}}
\newcommand{\cF}{\mathcal{F}}

\newcommand{\cM}{\mathcal{M}}
\newcommand{\cQ}{\mathcal{Q}}
\newcommand{\cR}{\mathcal{R}}
\newcommand{\cX}{\mathcal{X}}

\newcommand{\eps}{\varepsilon}

\newcommand{\p}{\partial}

\newcommand{\N}{\mathbb{N}}
\newcommand{\R}{\mathbb{R}}

\newcommand{\Z}{\mathbb{Z}}

\def\O{\Omega}

\newcommand{\abs}[1]{| #1 |}
\newcommand{\nor}[1]{\| #1 \|}

\DeclareMathOperator{\Div}{div}

\DeclareMathOperator{\argmin}{argmin}

\DeclareMathOperator{\supp}{supp}

\newcommand{\ptw}{\emph{Ptw}}
\newcommand{\rr}{\emph{RR}}
\newcommand{\fft}{\emph{FFT}}
\newcommand{\ffto}{\emph{FFTO}}
\newcommand{\fftten}{\emph{FFT10}}

\def\CC{{C\nolinebreak[4]\hspace{-.05em}\raisebox{.4ex}{\tiny\bf ++}}}

\title{ Error analysis of some nonlocal diffusion discretization schemes
\thanks{Supported by Spanish MCI Project MTM2017-87162-P.}}

\author{Gonzalo Galiano  \thanks{Dpt. of Mathematics, Universidad de Oviedo,
 c/ Calvo Sotelo, 33007-Oviedo, Spain ({\tt galiano@uniovi.es)}}}

\date{}

\begin{document}

\thispagestyle{plain}

\maketitle

\begin{abstract}
  We study three numerical approximations of solutions of nonlocal diffusion evolution problems which are inspired in algorithms for computing the bilateral denoising filtering of an image, and which are based on functional rearrangements and on the Fourier transform. Apart from the usual time-space discretization, these algorithms also use the discretization of the range of the solution (quantization). We show that the discrete approximations converge to the continuous solution in suitable functional spaces, and provide error estimates. Finally, we present some numerical experiments illustrating the performance of the algorithms, specially focusing in the execution time.   

\no\emph{Keywords: }Numerical approximation, nonlocal diffusion, functional rearrangement, Fourier transform, error analysis


\end{abstract}

\section{Introduction}

Let $T>0$ and $\O\subset\R^d$ $(d\geq 1)$ be an open and bounded set with Lipschitz continuous boundary, and consider the following problem:   Find $u:[0,T]\times\O\to\R_+$ such that 
\begin{align}
&  \p_t u(t,\bx)  =  \int_\O J(\bx,\by) A (u(t,\by)-u(t,\bx)) d\by +f(t,\bx,u(t,\bx)), \label{eq.eq}\\
&   u(0,\bx)=  u_0(\bx),  \label{eq.id}
\end{align}
 for $(t,\bx)\in Q_T=(0,T]\times \O$, and for some $u_0:\O\to\R$. 
The functions $J$ and $A$ are called, respectively, the \emph{spatial} and the \emph{range} kernels. 

Nonlocal diffusion problems of the type \fer{eq.eq}-\fer{eq.id} have many applications in modelling diffusion processes in which the interactions among particles occur at a distance, in contrast to local diffusion modelled by partial differential equations (PDE). 

They appear in a number of fields, among which image processing \cite{Buades2005, Gilboa2008} and population dynamics \cite{Hutson2003, Hetzer2011} are specially relevant. For the latter, the usual choice of the monotone range kernel $A(s)=\abs{s}^{p-2}s$, with $p\geq1$, gives rise to the so called \emph{nonlocal $p$--Laplacian} due to some similarities to the usual evolution $p$--Laplacian PDE. For the former, non-monotone but Lipschitz continuous range kernels like the Gaussian function are more frequently used. In both of them, the spatial kernel is assumed to be, at least, an integrable non-negative symmetric function.

In addition to the well-posedness theory, see the monograph by Andreu et al. \cite{Andreu2010} for the monotone case and \cite{Galiano2019, GalianoVelasco2019} for the Lipschitz continuous case, the numerical discretization of this type of problems has also been investigated, being most of the analyzed methods related to the usual techniques for PDE. For instance, P\'erez-Llano and Rossi \cite{Perez2011} analyze  a discrete spatial decomposition in terms of Lagrange basis while Tiang and Du \cite{Tiang2013} study, in addition, the standard finite element approach. In  \cite{Du2019}, Du et al.  consider a discontinuous Galerkin method. See the recent monographs by D'Elia et al. for a review on finite element based methods \cite{Delia2020B} and other methods \cite{Delia2020A} that also apply to non-integrable spatial kernels.

In this article we study some discretization techniques which were introduced in the image processing literature without performing a rigorous mathematical analysis. The purpose of these techniques is the fast computation of the approximations, being its accuracy relegated as of secondary importance. As stated by Weiss \cite{Weiss2006}, the fundamental property that concerns us is the runtime per pixel, as a function of the spatial kernel support. 

The first method we analyze is based on the work of Weiss \cite{Weiss2006}, later improved by Porikli \cite{Porikli2008} and generalized in \cite{Galiano2015A, Galiano2015B}. The main idea is changing the space of integration in \fer{eq.eq}. In the simplest discrete situation, the problem is to compute efficiently the filter
\begin{align}
\label{def:typebb}
Fu[\bi] = \sum_{\bj\in\bJ} J[\bi-\bj]A([u[\bj]-u[\bi]),	
\end{align}
where $\bJ\subset \Z^2$ is the set of pixels and $u:\bJ\to Q=\{0,1,\ldots,255\}$ is an intensity image. When $J$ is a \emph{box} kernel (taking a constant positive value inside a small square and zero outside) then the filter may be computed using the \emph{local histogram}, $h_\bi$, of $u$ in the box centered at each pixel $\bi\in\bJ$, 
\begin{align}
\label{def:irr}
Fu[\bi] = \sum_{k \in Q}  h_\bi(k) A(k-u[\bi]),
\end{align}
that is, changing the space of integration from the pixel space to the image value space. In the continuous case, and for general integrable spatial kernels, this change involves in a natural way the coarea formula, leading to the consideration of some functional rearrangements involving the image (decreasing rearrangement) and the spatial kernel (relative rearrangement). 
 
 The second method was introduced by Durand and Dorsey \cite{Durand2002} and later computationally improved by Yang et al. \cite{Yang2009, Yang2015}. Considering the collection of convolution filters
\begin{align}
\label{def:ifft}
F_k u[\bi] = \sum_{\bj\in\bJ} J[\bi-\bj]A([u[\bj]-k),	\qtext{for }k\in Q,
\end{align}
 one observes that $F u[\bi] = F_k u[\bi]$ if $u[\bi]=k$. The convolution filters may be efficiently computed by means of the fast Fourier transform. Moreover, a natural and fast approximation to the exact filter may be defined by calculating the convolution filters only for a small subset of $Q$ and then defining the approximation through linear interpolation.  
 
For image processing tasks, specially for image denoising, filters involving expressions like \fer{def:typebb} are known as \emph{bilateral filters} \cite{Yaroslavsky1985, Smith1997, Tomasi1998, Barash2002, Elad2002, Buades2006, Gilboa2008} and  are mostly applied in a single discrete time step, although other tasks like image segmentation and saliency detection may require more steps \cite{Galiano2015A, Galiano2016, Galiano2020}. For other applications, these algorithms may result interesting when accuracy is secondary to execution time, e.g. for real time computations. 

The error analysis we perform in this article fills a gap, already pointed out by Durand and Dorsey \cite{Durand2002}, with the aim of establishing a rigorous mathematical setting for these algorithms, usually formulated in the discrete case and for concrete spatial and range kernels, mostly the box and the Gaussian kernels. With a well-posedness theory at hand \cite{Andreu2010, Galiano2019}, we extend the ideas leading to the filters \fer{def:irr} and \fer{def:ifft} to the continuous evolution case of the form \fer{eq.eq}-\fer{eq.id}, which apart of adding a reaction term, it also allows the consideration of  general spatial and range kernels. More concretely, 
adopting the assumptions for the existence of solutions results \cite{Andreu2010, Galiano2019}, we suppose \textbf{(H)}:
\begin{enumerate}
 \item The \emph{spatial kernel} $J\in  BV(\R^d\times\R^d)$ is symmetric and non-negative. Without loss of generality, and for normalization porposes, we also assume
 \begin{align*}
\int_{\O\times\O} J(\bx,\by)d \bx d\by = 1.
 \end{align*}

\item The \emph{range kernel} $A \in W_{loc}^{1,\infty}(\R)$ satisfies  
\begin{align*}
& A(-s) = - A(s),\quad  A(s)s\geq0 \qtext{for all }s\in\R.
\end{align*}

\item The reaction function $f \in \big(L^\infty(Q_T)\cap BV(Q_T)\big) \times W_{loc}^{1,\infty}(\R)$ satisfies 
\begin{align*}
& \abs{f(\cdot,\cdot,s)}\leq C_f (1+\abs{s}) \text{ in }Q_T, \text{ for }s\in\R, \text{ and for some constant } C_f>0. 
\end{align*}

\item The initial datum $u_0\in L^\infty(\O)\cap BV(\O)$. 
\end{enumerate}  
Under the set of assumptions (H), the existence and uniqueness of a solution, $u$,  of \fer{eq.eq}-\fer{eq.id} such that
$u\in W^{1,\infty}(0,T;L^\infty(\O))\cap C([0,T];L^\infty(\O)\cap BV(\O))$,
is proven in \cite[Theorem 1]{Galiano2019}.  In particular,  it is shown that 
\begin{align}
\label{existence:bound}
 \nor{u}_{L^\infty(Q_T)} \leq C_\peq{M},
\end{align}
where $C_\peq{M}$ only depends on $T$, $\nor{u_0}_{L^\infty(\O)}$ and $C_f$. In addition, if 
\begin{align}
\label{ass.pos}
u_0\geq0 \qtext{in } \O \qtextq{and }f(\cdot,\cdot,0)\geq 0 \qtext{in }Q_T	
\end{align}
 then the solution is non-negative in $Q_T$. In the rest of the paper, we assume \fer{ass.pos}, for simplicity.

Let us make some observations. On one hand, 
the well-posedness of problem \fer{eq.eq}-\fer{eq.id} may be also guaranteed under other complementary conditions on the range kernel that replace (H)$_2$: $A\in C^{0,\alpha}_{loc}(\R)$ is a non-decreasing function, see \cite[Theorem 2]{Galiano2019}. In fact, the H\"older continuity may be weakened under some circumstances, see \cite{Andreu2010}. A numerical analysis of this case may be found in \cite{Perez2011}. Besides, we may consider more general conditions than those stated in (H). For instance, we may include a time-space dependence on the range kernel, or  a nonlocal reaction term in place or in addition to the local reaction term, $f$. With the aim of keeping some notational simplicity we have opted to deal with a simpler model, being the techniques described in this article easy to extend to the general case.

On the other hand, if the spatial regularity of $J$, $f$ and $u_0$ is improved, for instance to a $W^{1,p}$ Sobolev space, with $p\in(1,\infty]$,  then the corresponding spatial regularity of the solution is improved to the same space. Notice that one of the most remarkable properties of nonlocal diffusion evolution problems, when compared to their local diffusion versions, is the lack of a regularizing effect on the solutions. This is, the solutions are not more regular than the data.
	
Finally, observe that the important limit case of the nonlocal $1$--Laplacian equation is not covered by the assumptions (H). However, a proof of its well-posedness may be found in \cite{Andreu2010}. In contrast, (H) allows to handle the case $A(s)=\abs{s}^{p(s)-2}s$, this is, the nonlocal version of the problem $\p_t u - \Div\big({\abs{\nabla u}^{p(\nabla u)-2}\nabla u}\big) = f$, already employed in the literature \cite{Blomgren1997}, but for which the well-posedness is an open problem.

The contents of the rest of the article is the following. In Section~\ref{sec:general} we introduce the general discretization scheme, based on an explicit time Euler method and a piecewise constant approximation in time and  space. Since the only difference among the methods is in the quadrature of the integral term in \fer{eq.eq}, we give in this section the error estimates for the other terms. In Section~\ref{sec:methods} we provide the analysis of the three quadrature methods considered in the article: the pointwise method (\ptw), the method based in functional rearrangements (\rr) and  the method based on the Fourier transform (\fft). Finally, in Section~\ref{sec:numerics} we show the results of some numerical experiments illustrating the performance of each method for different choices of the space and range kernels and for a variety of discretization parameters. 

\section{A general approximation scheme}\label{sec:general}

The proof of the existence of solutions in \cite{Galiano2019} suggests the use of an explicit scheme to discretize \fer{eq.eq}-\fer{eq.id} due to the absence of a stability constraint for the time step in terms of the spatial mesh size, see Remark~\ref{rem:stab}. 

Thus, all the fully discrete schemes considered in this article are based on the following semi-discrete scheme for problem \fer{eq.eq}-\fer{eq.id}: Set  $\tau = T/N$, for some $N\in\N$ and $t_n= \tau n$ for $n=0,1,\ldots,N$. Then, for $\bx\in \O$, and $n=0,1,\ldots,N-1$, we set $u^{0}(\bx)=u_0(\bx)$ and 
\begin{align}
 &   u^{n+1}(\bx)=  u^n(\bx)+\tau \big( \cA(u^n)(\bx) + f_n(\bx, u^n(\bx))\big),\label{eq:evoln}
 \end{align}
where $f_n$ is some approximation of $f(t,\cdot,\cdot)$ for $t\in [t_n, t_{n+1})$ and where, for $v\in L^\infty(\O)$ and $\bx\in\O$, we define 
\begin{align}
\cA(v)(\bx) = \int_\O J(\bx,\by) A(v(\by)-v(\bx)) d\by. \label{def:operator}
\end{align}

We take $\O\subset\R^d$ as the hyper-cube $\O = (0,L)^d$, with $L\in\R$ and consider the uniform mesh of $\O$ of size $h$, $\cM_h(\O) = \{\bx_\bj\}_{\bj\in\bJ}$,  given by the centroids
\begin{align}
\label{space:disc}
\bx_\bj = \big(\bj - \frac{{\bf 1}}{2} \big)h ,\qtext{for } \bj \in \bJ = \{( j_1 ,\ldots, j_d ),\quad j_i=1,\ldots,J_\O\}, 
\end{align}
where ${\bf 1}\in \R^d$ is given as ${\bf 1}=(1,\ldots,1)$ and $J_\O=L/h$ is assumed to be integer.
Thus, the collection of hypercubes of volume $\abs{\O_\bj} = h^d$
centered at $\bx_\bj$
\begin{align*}
\O_\bj=\{\bx\in\O: \max_{k=1,\ldots,d}\abs{\bx_k-(\bx_\bj)_k}<h/2\} 
\end{align*}
induces a partition of $\O$.

Let $v_h$ be a generic (G) piecewise constant approximation of $v\in L^1(\O)$ defined by
\begin{align}
\label{def:vG}
 v_h^\peq{G}(\bx) = v_h^\peq{G}[\bj] \qtext{if }\bx\in\O_\bj, 
\end{align}
where, abusing on notation, we use the same symbol, $v_h^\peq{G}$, for the function $v_h^\peq{G}:\O\to\R$, whose evaluation at $\bx\in\O$ is denoted by $v_h^\peq{G}(\bx)$, and for the constant value this function takes in $\O_\bj$, which is denoted by $v_h^\peq{G}[\bj]$. The value $v_h^\peq{G}[\bj]$ will be specified for each discretization method.
Similarly, assume the following generic forms for the approximations of $\cA$ and $f$:
\begin{align}
\label{def:AfG}
\cA_h^\peq{G} (v_h^\peq{G})(\bx) = \cA_h^\peq{G} (v_h^\peq{G})[\bj] \qtext{if }\bx\in\O_\bj, \qtextq{and} f_{\tau h}^\peq{G}(t,\bx,s) = f_{n\bj}^\peq{G}(s) 
\end{align}
if $(t,\bx)\in [t_n,t_{n+1})\times\O_\bj$.
Then, we consider the discrete scheme given by,  for $\bj\in\bJ$ and $n=0,\ldots,N-1$,  
\begin{align}
& u^0[\bj] = u_{0h}^\peq{G}[\bj], \nonumber \\
& u^{n+1}[\bj]=  u^n[\bj]+\tau \big( \cA_h^\peq{G} (u_{\tau h}(t_n,\cdot))[\bj]
+ f_{n\bj}^\peq{G}(u^n[\bj]) \big),
\label{eq:scheme}
\end{align} 
where $u_{\tau h}$  is the piecewise constant  interpolator of the collection $\{u^n[\bj]\}_{\bj\in\bJ}$, defined as 
\begin{align*}
u_{\tau h}(t,\bx)=u^{n}[\bj] \qtext{if }(t,\bx)\in [t_n,t_{n+1})\times \O_\bj. 
\end{align*}
We also introduce here, for future reference, the corresponding
time piecewise linear interpolator
\begin{align*}
 \overline{u}_{\tau h}(t,\bx)=u^{n+1}[\bj] +\frac{t_{j+1}-t}{\tau}(u^n[\bj])-u^{n+1}[\bj]) \qtext{if }(t,\bx)\in [t_n,t_{n+1})\times \O_\bj, 
\end{align*}
and fix the initial datum approximations as $ u_{\tau h}(0,\cdot)= \overline{u}_{\tau h}(0,\cdot)= u_{0h}^\peq{G}$ in $\O$.

Observe that, for notational consistency, we should write $u_{\tau h}^\peq{G}$ and $u^{n\peq{G}}$ instead of $u_{\tau h}$ and $u^{n}$. For the sake of clearness, we drop this index from the approximating functions in the hope that the context will avoid any confusion.  

\subsection{Error bound for the common terms}

We shall deduce, for each discretization method, error bounds in the usual norm of $L^\infty(0,T;L^1(\O))$. The choice of this norm is motivated by the functional space where generalized differentials of $BV$ functions may be approximated by piecewise constant functions, see Lemma~\ref{lema:approx}. 

Subtracting the equations satisfied by $u$ and $\overline{u}_{\tau h}$ we get, for $(t,\bx)\in [t_n,t_{n+1})\times \O_\bj$, 
\begin{align*}
\p_t\big(u(t,\bx) -  \overline{u}_{\tau h}(t,\bx)\big) = & ~\cA (u(t,\cdot))(\bx) - \cA_h^\peq{G}(u_{\tau h}(t,\cdot))(\bx) \nonumber \\
 & + f(t,\bx,u(t,\bx)) -  f_{\tau h}^\peq{G}(t,\bx,u_{\tau h}(t,\bx)). 
\end{align*}
Integrating in $Q_t$ and using the triangle inequality,  we obtain the estimate
\begin{align}
\int_\O \abs{u(t,\bx)  -  u_{\tau h}(t,\bx)} d\bx \leq & ~ \int_\O \abs{u_0(\bx)  -  u_{0h}^\peq{G}(\bx)} d\bx+ 
\int_\O \abs{\overline{u}_{\tau h}(t,\bx) -  u_{\tau h}(t,\bx)} d\bx  \nonumber  \\
&+\int_{Q_t} \left|\cA (u(s,\cdot))(\bx)-\cA (u_{\tau h}(s,\cdot))(\bx)\right| d\bx ds \nonumber  \\ 
&+\int_{Q_t} \left|\cA  (u_{\tau h}(s,\cdot))(\bx)-\cA_h^\peq{G} (u_{\tau h}(s,\cdot))(\bx)\right| d\bx ds \nonumber  \\ 
& + \int_{Q_t} \abs{f(s,\bx,u(s,\bx)) -  f(s,\bx, u_{\tau h}(s,\bx))}d\bx ds  \nonumber \\
& + \int_{Q_t} \abs{f(s,\bx, u_{\tau h}(s,\bx)) -  f_{\tau h}^\peq{G}(s,\bx, u_{\tau h}(s,\bx))}d\bx ds  \nonumber \\
& =: I_1+\cdots+I_6. \label{conv:fv1}
\end{align}
From the assumptions (H) we deduce
\begin{align}
& \max\{I_3,I_5\} \leq C \int_{Q_t} \abs{u(s,\bx) -  u_{\tau h}(s,\bx)} d\bx ds. \label{est:f1}
\end{align}
The remaining terms of \fer{conv:fv1} will be estimated in the following section, according to the type of approximation considered.

\section{Discretization methods}\label{sec:methods}
\subsection{The Pointwise method}

\no\textbf{Definition of the approximation. }We start fixing the values of the approximations defined in \fer{def:vG} and \fer{def:AfG}. We set $G=P$ (\emph{Pointwise}), and 
\begin{align}
 &v_h^\peq{P}[\bj]=\frac{1}{h^d}\int_{\O_\bj} v(\bx)d\bx,\qtext{for }\bj\in\bJ, \label{def:trozos}\\
&f_{n\bj}^\peq{P}(s) = \frac{1}{\tau h^d} \int_{t_n}^{t_{n+1}}\int_{\O_\bj}f(t,\bx,s) d\bx dt, \qtext{for }n=0,\ldots,N,~\bj\in\bJ,~s\in\R,\label{def:discf}\\
& \cA_h^\peq{P}(v_h^\peq{P})[\bj] = h^d\sum_{\bk\in\bJ} J_h^\peq{P} [\bj,\bk] A(v_h^\peq{P}[\bk]-v_h^\peq{P}[\bj]), \qtext{for }\bj\in\bJ, \nonumber
\end{align}
with 
\begin{align}
\label{def:Jh}
 J_h^\peq{P}[\bj,\bk] = \frac{1}{ h^{2d}} \int_{\O_\bj}\int_{\O_\bk}J(\bx,\by) d\by d\bx, \qtext{for }\bj,\bk\in\bJ. 
\end{align}
Observe that  the definitions of $\cA_h^\peq{P}$ and $J_h^\peq{P}$ imply 
\begin{align*}
 \cA_h^\peq{P}(v_h^\peq{P})[\bj] = \frac{1}{h^d}\int_{\O_\bj} \cA(v_h^\peq{P})(\bx) d\bx, \qtext{if }\bx\in\O_\bj.
\end{align*}
The following lemma allows to estimate the approximation error between a $BV$ function and its piecewise constant approximation defined in terms of its mean values in $\O_\bj$, for $\bj\in\bJ$.
\begin{lemma}
\label{lema:approx}
Let $v\in BV(\O)$ and consider the function $v_h^\peq{P}$ defined by \fer{def:vG} and \fer{def:trozos}. Then there exists a constant $C>0$, independent of $v$ and $h$, such that	
\begin{align}
\nor{v_h^\peq{P}-v}_{L^1(\O)} \leq Ch  \abs{Dv}(\O), \label{est:PW}
\end{align}
where $\abs{Dv}(\O)$ is the variation of $v$ in $\O$. Similarly, we have 
\begin{align}
& \nor{J_h^\peq{P} - J}_{L^1(\O\times\O)} \leq C h \abs{DJ}(\O\times\O), \label{est:PW3}\\
&\nor{f_{\tau h}^\peq{P}(\cdot,\cdot,s) - f(\cdot,\cdot,s)}_{L^1(Q_T)} \leq C (\tau+h) \abs{Df(\cdot,\cdot,s)}(Q_T),\label{est:PW2}
\end{align}
where $J_h^\peq{P}(\bx,\by) = J_h^\peq{P}[\bj,\bk]$ if $\bx\in\O_\bj$ and $\by\in\O_\bk$.
\end{lemma}
\no\textbf{Proof. } 
According to Theorem~3.9 of \cite{Ambrosio2000}, if $v\in BV(\O)$ then there exists a sequence $\{v^\eps\}\subset C^\infty(\O)$ such that $v^\eps\to v$ in $L^1(\O)$ and $\nor{\nabla v^\eps}_{L^1(\O)}\to \abs{Dv}(\O)$ as $\eps\to0$. Let $v^\eps_h:\O\to\R$ be defined by 
\begin{align*}
v^\eps_h (\bx) = \frac{1}{h^d}\int_{\O_\bj} v^\eps(\bx)d\bx \qtextq{if} \bx\in\O_\bj. 	
\end{align*}
On one hand, we have
\begin{align*}
\int_\O\abs{v^\eps_h(\bx) -v_h^\peq{P}(\bx)}	d\bx & = \sum_{\bj\in\bJ}\int_{\O_\bj} \frac{1}{h^d}\left|\int_{\O_\bj} (v^\eps(\bx)-v(\bx))d\bx \right|\\
& \leq \int_\O\abs{v^\eps(\bx)-v(\bx)}d\bx\to 0
\end{align*}
as $\eps\to0$. On the other hand, since $v^\eps_h|_{\O_\bj}$ is the mean value of $v^\eps$ in $\O_\bj$, Poincar\'e's inequality in the square $\O_\bj$ of diameter $h$ gives  $\nor{v_h^\eps-v^\eps}_{L^1(\O_\bj)} \leq Ch  \nor{\nabla v^\eps}_{L^1(\O_\bj)}$, with $C$ independent of $h$ and $v^\eps$. Summing for $\bj\in\bJ$, we deduce $\nor{v_h^\eps-v^\eps}_{L^1(\O)} \leq Ch  \nor{\nabla v^\eps}_{L^1(\O)}$. Thus, we obtain
\begin{align*}
\nor{v_h^\peq{P}-v}_{L^1(\O)} \leq  \nor{v_h^\peq{P}-v^\eps_h}_{L^1(\O)} +\nor{v^\eps-v}_{L^1(\O)} + Ch  \nor{\nabla v^\eps}_{L^1(\O)},
\end{align*}
and we deduce the result by taking the limit $\eps\to0$. The estimates \fer{est:PW3} and \fer{est:PW2}  are proven in a similar manner, replacing $\O$ by the corresponding product space. $\Box$

Here, and in what follows, $C$ denotes a positive constant  which may change of value among different expressions, but that is always independent of $\tau$, $h$ and other approximation parameters ($q$ and $r$, see \fer{rangedisc}).

The following lemma states that, as well as the solution $u$ of problem \fer{eq.eq}-\fer{eq.id}, the sequence of approximating functions $u_{\tau h}$ remains uniformly bounded in $L^\infty(Q_T)$. 
The result is a consequence of the Lipschitz continuity of $A$ and the boundedness of $u$, which imply the sublinearity of $A$ in $[0,\nor{u}_{L^\infty(Q_T)}]$.

The proof may be easily adapted to the discretizations introduced in Sections 	\ref{sec:rr} and \ref{sec:fourier}. We, therefore, omit the corresponding proofs.
\begin{lemma}
\label{lemma:linfty}
 There exists a constant $C>0$, independent of $\tau$ and $h$, such that 
 \begin{align}
\label{est:uniform}
\nor{u_{\tau h}}_{L^\infty(Q_T)} \leq C.
\end{align}
\end{lemma}

\no\textbf{Proof. }
Let $u$ be the unique solution of problem \fer{eq.eq}-\fer{eq.id}. Since it satisfies  the bound \fer{existence:bound}, we deduce from (H)$_2$ the existence of a positive constant $C_A$ such that $\abs{A(s)} \leq C_As$ for $s\in [0,C_\peq{M}]$. We may assume that this property is satisfied for all $s\in\R$ since, on the contrary, redefining $A$ in $s\in\R\backslash[0,C_\peq{M}]$ as $A(s)=C_A s$, the solution to problem \fer{eq.eq}-\fer{eq.id} is still the same function than for the original range kernel $A$. Hence, from now on we assume, without loss of generality, 
\begin{align*}
\abs{A(s)}\leq C_A \abs{s} \qtextq{for all}s\in\R.	
\end{align*}
Consider the change of unknown $w^n=e^{-\mu t_n}u^n$, for some $\mu>0$ to be fixed. Then, $w^n$ satisfies
\begin{align*}
& w^0[\bj] = u^0[\bj], \nonumber \\
& w^{n+1}[\bj]=  e^{-\mu\tau}w^n[\bj]+\tau e^{-\mu t_{n+1}} \big( \cA_h^\peq{P} (e^{\mu t_{n}}w_{\tau h}(t_n,\cdot))[\bj]
+ f_{n\bj}^\peq{P}(e^{\mu t_{n}}w^n[\bj]) \big),
\end{align*} 
for $\bj\in\bJ$. We have 
\begin{align*}
	 & \abs{\cA_h^\peq{P} (e^{\mu t_{n}}w_{\tau h}(t_n,\cdot))[\bj]}\leq e^{\mu t_{n}}2C_A\max_{\bk\in\bJ} \abs{w^n[\bk]} h^d\sum_{\bk\in\bJ} J_h^\peq{P} [\bj,\bk] ,\\
	 & \abs{f_{n\bj}^\peq{P}(e^{\mu t_{n}}w^n[\bj])} \leq e^{\mu t_{n}} C_f (1+\max_{\bk\in\bJ} \abs{w^n[\bk]} ),
\end{align*}
implying
\begin{align*}
& \abs{w^{n+1}[\bj]}\leq  e^{-\mu\tau} \Big(\tau C_f+\big(1 +\tau(2C_A+C_f)\big) \max_{\bk\in\bJ} w^n[\bk]\Big).
\end{align*} 
Therefore, taking $\mu>2C_A+C_f$ and returning to the original unknown this differences inequality yields the uniform estimate $\max_{\bk\in\bJ}  \abs{u^n[\bk]} \leq C$, with $C$ depending only on $\max_{\bk\in\bJ}  \abs{u^0[\bk]} $, $T$ and the constants $C_A$ and $C_f$. $\Box$

\begin{remark}
\label{rem:stab}
Since $A(s)$ and $f(\cdot,\cdot,s)$ are Lipschitz continuous in the interval
 $I=(0,C)$, with $C$ given in \fer{est:uniform}, a sufficient condition for the stability of the numerical discretization \fer{eq:scheme} is given by 
\begin{align}
\label{cond:est}
\tau < \frac{1}{2\nor{A}_{W^{1,\infty}(I)}+\nor{f}_{L^\infty(Q_T)\times W^{1,\infty}(I)}},	
\end{align}
which is similar to that found by P\'erez-Llano and Rossi \cite[Proposition 2.9]{Perez2011} when $A(s)=\abs{s}^{p-2}s$ and $f=0$. In fact, \fer{cond:est} ensures that the comparison principle is satisfied in the discrete setting of that problem. Notice that, in contrast to the corresponding local diffusion problem, \fer{cond:est} does not depend on the spatial mesh size.
\end{remark}

\no\textbf{Remaining estimates. }Estimate \fer{est:PW} applied to the initial datum yields
\begin{align}
\label{est:id}
 I_1 \leq Ch.
\end{align}
By direct computation  we get, for $(t,\bx)\in [t_n,t_{n+1})\times \O_\bj$, 
\begin{align*}
\overline{u}_{\tau h}(t,\bx) - u_{\tau h}(t,\bx) & = 
 (t-t_n)\frac{u^{n+1}[\bj]-u^n[\bj]}{\tau}  \\
& = (t-t_n)\big(\cA_h^\peq{P}(u_{\tau h}(t,\cdot))(\bx)
 + f_{\tau h}^\peq{P}(t,\bx, u_{\tau h}(t,\bx)) \big).
\end{align*}
Using \fer{est:uniform} and that $\abs{t-t_n}<\tau$ for $t\in[t_n,t_{n+1})$, we deduce 
\begin{align}
\label{est:time}
I_2  \leq C\tau .
\end{align}
For $(s,\bx)\in[t_n,t_{n+1})\times\O_\bj$, we have 
\begin{align*}
\cA(u_{\tau h}(s,\cdot))(\bx) & - \cA_h^\peq{P}(u_{\tau h}(s,\cdot))(\bx) \\
& = \sum_{\bk\in\bJ} A(u_n[\bk]-u_n[\bj] ) \Big(\int_{\O_\bk} J(\bx,\by) d\by - h^d J_h^\peq{P}[\bj,\bk] \Big) \nonumber\\
& = \sum_{\bk\in\bJ} A(u_n[\bk]-u_n[\bj] ) \Big(\int_{\O_\bk}  \big(J(\bx,\by)  -  J_h^\peq{P}(\bx,\by)\big) d\by \Big). 
\end{align*}
From \fer{est:PW3} and \fer{est:uniform}, we obtain
\begin{align}
  I_{4}  \leq C   \int_{\O} \int_{\O} \big|J(\bx,\by) d\by - J_h^\peq{P}(\bx,\by) \big| d\by d\bx \leq Ch^2, \label{est:A2}
\end{align}
and, in view of \fer{est:PW2}, we deduce
\begin{align}
\label{est:f2}
  I_6  \leq C\tau h.
\end{align}
Using estimates \fer{est:f1}, \fer{est:id}, \fer{est:time},  \fer{est:A2} and \fer{est:f2} in \fer{conv:fv1}, we finally obtain 
\begin{align*}
\int_\O \abs{u(t,\bx)  -  u_{\tau h}(t,\bx)} d\bx \leq &  ~ C(\tau + h ( \tau + h))    \\ &+ C\int_{Q_t} \abs{u(s,\bx)  -  u_{\tau h}(s,\bx)} d\bx ds,
\end{align*}
and, for $\tau,h<1$,  Gronwall's lemma implies 
\begin{align}
 \nor{u-u_{\tau h}}_{L^\infty(0,T;L^1(\O))} \leq C(\tau + h). \label{eb.ptw}
\end{align}

\subsection{Using functional rearrangements}\label{sec:rr}

This first non-trivial discretization of the nonlocal functional \fer{def:operator} makes use of functional rearrangements. We give a short and formal description of the deduction of the main formulas we use, and refer the reader to \cite{Galiano2015B} and its references for the details.

Under suitable regularity assumptions on the functions $v$ and $g$ defined in $\O$, the coarea formula states 
\begin{equation}
\label{coarea}
 \int_\O g(\by)\abs{\nabla v(\by)} d\by =\int_{-\infty}^{\infty} \int_{v=s} g(\by) d\Gamma(\by) ds.
\end{equation}
Taking $ g(\by)= J(\bx,\by) A(v(\by)-v(\bx))  /\abs{\nabla v(\by)}$,
and assuming  $v(\bx) \in [0,Q]$ for all $\bx\in\O$ we get from \fer{def:operator} and \fer{coarea},
\begin{align}
 \cA(v)(\bx) =\int_{0}^{Q}  A(s-v(\bx))   
  \int_{v=s}  \frac{J(\bx,\by)}{\abs{\nabla v(\by)}} d\Gamma(\by) ds. \label{F.F}
\end{align}
In this expression, the inner integral is computed relative to the level sets of $v$, but independently of the \emph{ordering} or such sets. This is the main motivation to introduce a reformulation of \fer{def:operator} in terms of the decreasing and the relative rearrangements.
\subsection*{The decreasing and the relative rearrangements}

Let us denote by $\abs{E}$ the Lebesgue measure of any measurable set $E$.
For a Lebesgue measurable function $u:\O\to\R$, the function 
$q\in\R\to m_v(q) = \abs{\{\bx \in\O : v(\bx) >q\}}$ is called the \emph{distribution function} corresponding to $v$. 

Function $m_v$ is non-increasing and therefore admits a unique  generalized inverse, called the
\emph{decreasing rearrangement}. This inverse takes the usual pointwise meaning when 
the function $v$ has no flat regions, i.e. when $\abs{\{\bx \in\Omega : v(\bx) =q\}} =0$ for any $q\in\R$. In general, 
the decreasing rearrangement $v_*:[0,\abs{\O}]\to\R$ is given by:
\begin{equation*}
v_*(s) =\left\{
\begin{array}{ll}
 {\rm ess}\sup \{v(\bx): \bx \in \O \} & \qtext{if }s=0,\\
 \inf \{q \in \R : m_v(q) \leq s \}& \qtext{if } s\in (0,\abs{\O}),\\
 {\rm ess}\inf \{v(\bx): \bx \in \O \} & \qtext{if }s=\abs{\O}.
\end{array}\right.
 \end{equation*}
We shall also use the notation $\O_*=(0,\abs{\O})$. 

The notion of rearrangement of a function was introduced by Hardy, Littlewood and Polya \cite{Hardy1964}. 
Its most remarkable property is the  equi-measurability:  for $v\in L^1(\O)$ and for any Borel function $f:\R\to\R_+$
\begin{equation}
\label{prop.1}
 \int_\O f(v(\by))d\by = \int_{\O_*} f(v_* (s))ds.
\end{equation}

Returning to formula \fer{F.F} and assuming that $v_*$ is strictly decreasing, we may introduce the change of variable $t=v_{*}(s)$ and use \fer{prop.1} to obtain
\begin{align}
 \cA(v)(\bx) & =-\int_{\O_*}  A(v_*(s)-v(\bx)) \frac{d v_* (s)}{ds}  \int_{v=v_*(s)}  \frac{J(\bx,\by)}{\abs{\nabla v(\by)}} d\Gamma(\by) ds  \nonumber \\
 & =\int_{\O_*} J(\bx,\cdot)_{*v} (s) A(v_*(s)-v(\bx))  ds.
 \label{chco}
\end{align}
Here, the notation $\phi_{*v}$ stands for the \emph{relative rearrangement of $\phi$ with respect to $v$} which, 
under regularity conditions, may be expressed as \cite{Rakotoson2008}
\begin{equation*}
 \phi_{*v}(s)  = -\frac{d v_* (s)}{ds} \displaystyle\int_{v=v_*(s)}  \frac{\phi(\by)}{\abs{\nabla v(\by)}} d\Gamma(\by).
\end{equation*}
More in general,   the relative rearrangement is understood as the weak $L^p(\O_*)$ directional derivative 
(weak* $L^\infty(\O_*)$, if $p=\infty$),
\begin{equation*}
 \phi_{*v}= \lim_{t\to0}\frac{(v+t\phi)_* - v_*}{t}. 
\end{equation*}
The relative rearrangement was introduced by Mossino and Temam \cite{Mossino1981}. We refer the reader to the  monograph written by Rakotoson \cite{Rakotoson2008} as a fundamental reference in this field.

In  \cite[Theorem 1]{Galiano2015B},  the equivalence between $\cA(v)(\bx)$ given by \fer{def:operator} and 
the expression \fer{chco} was rigorously established. More concretely, assuming (H) and $v\in L^1(\O)$, and defining 
 \begin{align}
\label{def.ubar}
 \cA_*(v)(\bx)  = \int_{\O_*}  J(\bx,\cdot)_{*v} (s) A(v_*(s)-v(\bx))  ds,
\end{align}
we have $\cA_*(v)(\bx) = \cA(v)(\bx)$ for a.e. $\bx \in \O$.
The main advantage of this formulation is that the multi-dimensional set of integration in \fer{def:operator} is reduced to
a one-dimensional set in \fer{def.ubar}. However, we have to pay this reduction by computing the 
relative rearrangement $J(\bx,\cdot)_{*v}$ for a.e. $\bx\in\O$. As we shall see, the price to pay is not too high when the functions involved in the computations are piecewise constant.

\newpage

\no\textbf{Definition of the approximation. }In this section, we assume $J\in L^\infty(\O\times\O)$.  We set $G=R$ (\emph{Rearrangement}). 
The definitions of $u_{0h}^\peq{R}$ and $f_{n\bj}^\peq{R}$ are the same as in \fer{def:trozos} and \fer{def:discf}.
For those of $\cA_h^\peq{R}$ and $J_h^\peq{R}$ we need to introduce some additional tools.

Let $L_\peq{Q}\in\N$ and consider the uniform mesh, $\{q_i\}_{i=0}^{L_\peq{Q}}$, of an interval $[C_m,C_\peq{M}] \supset[\inf(u),\sup(u)]$, where  $u$ is the solution of problem \fer{eq.eq}-\fer{eq.id}. We choose $q_i$  decreasingly ordered, that is, given by 
\begin{align}
\label{rangedisc}
q_0=C_\peq{M}, \quad  q_i=q_0-iq,\qtext{for }i=1,\ldots,L_\peq{Q},\qtext{with  } q= \frac{C_\peq{M}-C_m}{L_\peq{Q}}.
\end{align}
Similarly, let $L_\peq{R}\in\N$ and consider the uniform mesh, $\{r_i\}_{i=0}^{L_\peq{R}}$, of the interval $[\inf(J),\sup(J)]$. We also take $\{r_i\}_{i=0}^R$ decreasingly ordered, 
\begin{align*}
r_0=\sup(J), \quad  r_j=r_0-jr,\qtext{for }j=1,\ldots,L_\peq{R},\qtext{with  } 
 r= \frac{\sup(J)-\inf(J)}{L_\peq{R}}. 
\end{align*}
We introduce the decomposition of a real number, $\rho$, in its quantized and remainder parts relative to a mesh $\{s_i\}_{i=0}^{L_\peq{S}}$, for $S=Q,R$, given by  $\rho = \cQ_\peq{S}(\rho)+ \cR_\peq{S}(\rho)$, where
\begin{align}
\label{parteentera}
\cQ_\peq{S}(\rho) = s_m,\qtext{with } m = \underset{0\leq i \leq S}{\argmin} ~ \abs{\rho-s_i} , \qtext{and } \abs{\cR_\peq{S}(\rho)}\leq \frac{1}{2} s,
\end{align}
being $s$ the size of the mesh.

Consider the functional 
\begin{align*}
  &\cA^\peq{R}(v)(\bx) = \int_\O \cQ_\peq{R}(J(\bx,\by))A(\cQ_\peq{Q}(v(\by))- \cQ_\peq{Q}(v(\bx))) d\by.
\end{align*}
We have 
\begin{align}
\cA (v)(\bx) = & \cA^\peq{R} (v)(\bx) + \eps(\bx), \qtext{with} \label{def:ARR}\\
  \eps(\bx) = & \int_\O \cQ_\peq{R}(J(\bx,\by))\Big( A(v(\by)-v(\bx)) - A(\cQ_\peq{Q}(v(\by))- \cQ_\peq{Q}(v(\bx))) \Big)d\by \nonumber\\
  & +  \int_\O \big( J(\bx,\by)- \cQ_\peq{R}(J(\bx,\by))\big) A(v(\by)-v(\bx)) d\by. \nonumber 
\end{align}
Since $J\in L^1(\O\times\O)$ and $A$ is Lipschitz continuous, if $v\in L^\infty(\O)$ we deduce
\begin{align}
& \abs{\eps(\bx)} \leq C(q+r). \label{Aq1:split2}
\end{align}
The quantized functions are simple finite functions, which may be expressed as 
\begin{align}
\label{def:trunc}
\cQ_\peq{Q}(v(\bx)) = \sum_{i=0}^{L_{\peq{Q}}} q_i \cX_{U_i}(\bx),\quad  \cQ_\peq{R}(J(\bx,\by)) = \sum_{m=0}^{L_\peq{R}} r_m \cX_{V_m}(\bx,\by),
\end{align}
where $U_i$ and $V_m$ are the level sets of $\cQ_\peq{Q}\circ v$ and $\cQ_\peq{R}\circ J$. 
Then, as shown in \cite[Theorem~1]{Galiano2015B}, we have that if $\bx\in U_k$, for some $k=0,\ldots,Q$, then
\begin{align}
\label{def.NLRD}
\cA^\peq{R}(v)(\bx)= \sum_{i=0}^{L_\peq{Q}} A(q_i-q_k) \sum_{m=0}^{L_\peq{R}}  r_m \mu_m^i(\bx),
\end{align}
where 
$\mu_m^i(\bx)=\text{meas}\{\by \in U_i : \cQ_\peq{R}(J(\bx,\by))=r_m \}$. 

If $v\equiv v_h$ is a piecewise constant function defined as in \fer{def:trozos}, then the level sets $U_i$ introduced in \fer{def:trunc} are unions of hypercubes $\O_\bj$.
Inspired in the expression \fer{def.NLRD}, we define the discrete nonlocal diffusion operator as, for $\bj \in \bJ$ such that $\O_\bj\subset U_k$, 
\begin{align}
\label{def:AGmu}
\cA_{h}^\peq{R} (v_h)[\bj] 
= \sum_{i=0}^{L_\peq{Q}} A(q_i-q_k) \sum_{m=0}^{L_\peq{R}}  r_m \mu_m^i[\bj],\qtextq{with}
\mu_m^i[\bj] = \frac{1}{h^d}\int_{\O_\bj}  \mu_m^i(\bx)d\bx.
\end{align}
Finally, $J_h^\peq{R}$ is defined like \fer{def:Jh}, with $J$ replaced by $\cQ_\peq{R}\circ J$.
\bigskip

\no\textbf{Remaining estimates. } We first notice that, abusing on notation, we still denote the approximation of the continuous solution as $u_{\tau h}$, even if, in fact, it also depends on the range quantization size meshes, $q$ and $r$. We proceed  as in Section~\ref{sec:fourier}.

The estimates for $I_1, I_2$ and $I_6$ are the same as in \fer{est:id}, \fer{est:time} and \fer{est:f2}. Thus, it only rests to estimate $I_4$. Using \fer{def:ARR} and \fer{Aq1:split2}, we get
\begin{align}
 \left|\cA (u_{\tau h}(s,\cdot))(\bx) -\cA_h^\peq{R} (u_{\tau h}(s,\cdot))(\bx)\right| \leq &
 \left|\cA^\peq{R} (u_{\tau h}(s,\cdot))(\bx) -\cA_h^\peq{R} (u_{\tau h}(s,\cdot))(\bx)\right| \nonumber\\& + C(q+r). \label{est:ARR4A}
\end{align}
Using \fer{def.NLRD} we get, for $(s,\bx)\in [t_n,t_{n+1})\times \O_\bj$ such that $\O_\bj\subset U_k^n=\{\by\in\O: u^n(\by) = q_k\}$,
\begin{align}
 \cA^\peq{R} (u_{\tau h}(s,\cdot))(\bx)& -\cA_h^\peq{R} (u_{\tau h}(s,\cdot))(\bx) =  \sum_{i=0}^{L_\peq{Q}} A(q_i-q_k)\sum_{m=0}^{L_\peq{R}}  r_m (\mu_m^{i,n}(\bx)-\mu_m^{i,n}[\bj]), \label{id:sect}
\end{align}
where $\mu_m^{i,n}$ is defined like $\mu_m^i$ with $U_i$ replaced by $U_i^n$, see \fer{def:AGmu}. We have
\begin{align*}
 \sum_{m=0}^{L_\peq{R}}  r_m \mu_m^{i,n}[\bj] = \frac{1}{h^d} \sum_{m=0}^{L_\peq{R}} r_m \int_{\O_\bj}  \mu_m^{i,n}(\bx)d\bx  = \frac{1}{h^d}  \int_{\O_\bj} M(\bx) d\bx,
\end{align*}
with $M(\bx)=\sum_{m=0}^{L_\peq{R}} r_m   \mu_m^{i,n}(\bx)$. 
Therefore, 
\begin{align*}
 \sum_{m=0}^{L_\peq{R}}  r_m (\mu_m^{i,n}(\bx)-\mu_m^{i,n}[\bj]) = 
 M(\bx)- \frac{1}{h^d}\int_{\O_\bj} M(\bx)d\bx.
\end{align*}
Let $\bJ_i^n = \{\bk\in\bJ : \O_\bk \subset U_i^n\}$ so that  $U_i^n = \displaystyle\cup_{\bk\in\bJ_i^n}\O_\bk$. Then 
\begin{align*}
 M(\bx) &  = \sum_{m=0}^{L_\peq{R}}  r_m  \sum_{\bk\in\bJ_i^n} 
 \text{meas}\{\by \in \O_\bk : \cQ_\peq{R}(J(\bx,\by))=r_m \}  \\
 & = \sum_{\bk\in\bJ_i^n} \int_{\O_\bk} \sum_{m=0}^{L_\peq{R}}  r_m ~
 1_{\{\by \in \O_\bk : \cQ_\peq{R}(J(\bx,\by))=r_m \}} (\by) d\by \\ 
 & = \sum_{\bk\in\bJ_i^n} \int_{\O_\bk}  \cQ_\peq{R}(J(\bx,\by)) d\by 
 = \int_{U_i^n} \cQ_\peq{R}(J(\bx,\by)) d\by .
\end{align*}
Since $J\in BV(\O\times\O)$, we deduce that $M\in BV(\O)$ and, therefore, from \fer{id:sect},  the boundedness of $A$ in $[0,\sup u_{\tau h}]$,  and the estimate \fer{est:PW} applied to $M$ we deduce
\begin{align*}
 \left|\cA^\peq{R} (u_{\tau h}(s,\cdot))(\bx) -\cA_h^\peq{R} (u_{\tau h}(s,\cdot))(\bx) \right| \leq Ch, 
\end{align*}
which, together with \fer{est:ARR4A} yields
\begin{align}
\label{est:A2RR}
 I_4 \leq C(h+q+r).
\end{align}
We conclude like in the previous section: using estimates \fer{est:f1}, \fer{est:id}, \fer{est:time},   \fer{est:f2} and \fer{est:A2RR} in \fer{conv:fv1}, we obtain, for $\tau,h<1$ ,
\begin{align*}
\int_\O \abs{u(t,\bx)  -  u_{\tau h}(t,\bx)} d\bx \leq & ~ C(\tau + h +q+r)   \\& + C\int_{Q_t} \abs{u(s,\bx)  -  u_{\tau h}(s,\bx)} d\bx ds.
\end{align*}
Gronwall's lemma implies 
\begin{align}
 \nor{u-u_{\tau h}}_{L^\infty(0,T;L^1(\O))} \leq C(\tau + h+q+r). \label{eb.rr}
\end{align}

\subsection{Using the Fourier transform}\label{sec:fourier}

 For notational simplicity, we relocate the set $\O = (0,L)^d$ to $\O_\peq{F} = (-L/2,L/2)^d$, where we introduce
the mesh $\cM_h(\O_\peq{F}) = \{\bx_\bj\}_{\bj\in\bJ_\peq{F}}$ similar to \fer{space:disc}, but with 
\begin{align}
 \bJ_\peq{F} = \left\{\bj=( j_1 ,\ldots, j_d ),\quad j_i=-\frac{J_\O}{2},\ldots,\frac{J_\O}{2},\quad i=1,\ldots,d \right\}, \label{mesh:fft}
\end{align}
with $J_\O$ even, for simplicity.
In this section, we assume that the spatial kernel is of convolution form, $J(\bx,\by) = w(\bx-\by)$, and replace Assumption $(H_1)$ by:
\begin{enumerate}
\item [${\bf (H)_F}$]  $w\in BV(\O_\peq{F})$ is even and non-negative.  In addition, either (i) w is periodic in $\O_\peq{F}$ and $w\in L^\infty(\O_\peq{F})$, or (ii) $w$ is  compactly supported in $\O_\peq{F}$.
\end{enumerate}
Observe that to have $w(\bx-\by)$ well defined for $\bx,\by\in\O_\peq{F}$, we need to be able to evaluate this function in $B=(-L,L)^d$. However, in $(H)_F$,  $w$ is only defined in $\O_\peq{F}$. We remove this obstacle by considering, according to the assumptions in $(H)_F$,  either (i) its periodic extension to $B$, or (ii) its extension by zero to $B$.
Correspondingly, we also extend the mesh and set 
$\cM_h(B) = \{\bx_\bj\}_{\bj\in\bJ_\peq{B}}$ similar to \fer{space:disc}, but with 
\begin{align*}
 \bJ_\peq{B} = \left\{\bj=( j_1 ,\ldots, j_d ),\quad j_i=-J_\O,\ldots,J_\O,\quad i=1,\ldots,d \right\}. 
\end{align*}
Thus,  
if $\bx\in\O_\bj$ and $\by\in\O_\bk$, with $\bj,\bk\in\bJ_\peq{F}$,  we have $\bx-\by \in \O_{\bj-\bk}$, with $\bj-\bk \in \bJ_\peq{B}$.

\bigskip

\no\textbf{Definition of the approximation. } We set $G=F$ (\emph{Fourier}). Consider the functional 
\begin{align*}
 &\cA^\peq{F} (v)(\bx) = \int_{\O_\peq{F}} w(\bx-\by)A(v(\by)- \cQ_\peq{Q}(v(\bx))) d\by, 
\end{align*}
with $\cQ_\peq{Q}$ introduced in \fer{parteentera} for the uniform mesh of the range of $u$ given in \fer{rangedisc}.
We have  
\begin{align}
\cA v(\bx) = & \cA^\peq{F} (v)(\bx) + \eps(\bx), \qtext{with} \label{def:AF}\\
  \eps(\bx) = & \int_{\O_\peq{F}} w(\bx-\by)\Big( A(v(\by)-v(\bx)) - A(v(\by)- \cQ_\peq{Q}(v(\bx))) \Big)d\by.  \nonumber 
\end{align}
Since $w\in L^1(\O_\peq{F})$ and $A$ is Lipschitz continuous, if $v\in L^\infty(\O_\peq{F})$ we deduce
\begin{align}
& \abs{\eps(\bx)} \leq Cq. \label{Aq1:split}
\end{align}
We also introduce the functional
\begin{align*}
& \cA^{\peq{F},i} (v)(\bx) = \int_{\O_\peq{F}} w(\bx-\by)A(v(\by)- q_i) d\by \qtext{if }\cQ_\peq{Q}(v(\bx))=q_i,
\end{align*}
for some $i=0,\ldots,Q$, 
which may be rewritten as the convolution
\begin{align}
\label{def.conv}
 \cA^{\peq{F},i} (v)(\bx)=
  w * H^i(v)(\bx),\qtext{with } H^i(v)(\by)=A(v(\by)-q_i).
\end{align}
If $v$ and $w$ are periodic functions in $\O_\peq{F}$ then $\cA^{\peq{F},i}$ may be computed  using the Fourier transform, this is,
\begin{align}
\label{def.convFo}
 \cA^{\peq{F},i} (v)=\cF^{-1}(\cF(w)\cF(H^i(v))).
\end{align}
 The natural discretization of \fer{def.convFo} is given in terms of the discrete Fourier transform, $\cF_d$. To employ this tool, the discrete sequence generated by the scheme \fer{eq:scheme}  must be periodic, property that we obtain if the piecewise constant approximations of $u_0$, $f$ and $w$ are  periodic. Since the spatial domain is an hyper-cube, we may do this by averaging the function values corresponding to opposed faces of the cube. For instance, if $\O\subset\R^2$, we split 
$\O_\peq{F}$ in three disjoint components: the \emph{interior} $\O_\peq{FI}^h=[-L/2+h,L/2-h]\times [-L/2+h,L/2-h]$,  the \emph{corners} 
\begin{align*}
	\O_\peq{FC}^h = \left\{\frac{L}{2}(z_1,z_2): z_j\in \{-1,1\} \text{ for }j=1,2\right\},
\end{align*}
 and the \emph{sides without the corners}, $\O_\peq{FS}^h = \O_\peq{F} \backslash (\O_\peq{FI}^h\cup\O_\peq{FC}^h)$. Then, we define the periodic constant-wise approximation of $u_0$ as
\begin{align*}
 u_{0h}^\peq{F}(\bx) = \begin{cases}
                    \dfrac{1}{h^2}\displaystyle\int_{\O_\bj}u_{0}(\bx) d\bx& \text{if } \bx \in  \O_\bj\cap \O_\peq{FI}^h \qtext{for some }\bj\in\bJ_\peq{F},\\
                    \dfrac{1}{2h^2}\displaystyle\int_{\O_\bj \cup \O_{\bj^s}}u_{0}(\bx)d\bx & \text{if } \bx \in  \O_\bj\cap \O_\peq{FS}^h \qtext{for some }\bj\in\bJ_\peq{F},\\
                    \\
                    \dfrac{1}{4h^2}\displaystyle\int_{\O_\peq{FC}^h}u_{0}(\bx)d\bx & \text{if } \bx \in  \O_\peq{FC}^h,
                   \end{cases}
\end{align*}
where $\bj^s$ is the node opposed to $\bj$. This is, if $\bj = (\pm1,j_2)$ then $\bj^s = (\mp 1,j_2)$, and if  $\bj = (j_1,\pm1)$ then $\bj^s = (j_1,\mp 1)$, for $j_1,j_2\neq -1,1$. 
We define in a similar way the periodic approximations of $f_{\tau h}^\peq{F}$ and $w_h^\peq{F}$.  The extension to higher dimensional domains is straightforward.

For the nonlocal diffusion operator,  we define, for  $\bj\in\bJ_\peq{F}$, 
\begin{align}
\label{def:cacfq}
 \cA_{h}^\peq{F} (v_h)[\bj] &=  \cA_{h}^{\peq{F},i} (v_h)[\bj] \qtextq{if} \cQ_\peq{Q}(v_h[\bj])=q_i, \\
& \qtextq{with} \cA_{h}^{\peq{F},i} (v_h)[\bj]=\cF_d^{-1}(\cF_d(w_h^\peq{F})\cF_d(H^i(v_h))) [\bj].\nonumber
\end{align}
Due to the periodicity imposed on the piecewise approximations of the data, we may use the circular convolution theorem to obtain
\begin{align*}
 \cF_d^{-1}(\cF_d(w_h^\peq{F})\cF_d(H^i(v_h))) [\bj] =  w_h^\peq{F} \circledast H^i(v_h) [\bj] ,
 \end{align*}
yielding a discrete scheme which mimics the convolutional form of the continuous problem whenever the data $u_0$, $w$ and $f$ are periodic functions in $\O_\peq{F}$, see \fer{def.convFo}.  However, notice that such  periodicity has not been  assumed and that the identity \fer{def.convFo} is, in general, not valid.

\bigskip

\no\textbf{Remaining estimates. } 
The estimate for $I_6$ is the same than  \fer{est:f2}. Thus, it only rests to estimate $I_1$ and $I_3$. 
We have 
\begin{align*}
 I_1= \int_{\O_\peq{FI}^h} \abs{u_0-u_{0h}^\peq{F}} + \int_{\O\backslash \O_\peq{FI}^h} \abs{u_0-u_{0h}^\peq{F}}.  
\end{align*}
Since the restriction of a $BV$ function to a smooth set is still a $BV$ function, we may use \fer{est:PW} to estimate the first term of the right hand side by $Ch$.
For the second (the \emph{border}), we have 
\begin{align}
\label{est:u0FB}
\int_{\O_\peq{F}\backslash \O_\peq{FI}^h} \abs{u_0-u_{0h}^\peq{F}} \leq C \nor{u_0}_{L^\infty(\O_\peq{F})} \text{Per}(\O_\peq{F}) h  ,
\end{align}
implying that 
\begin{align}
\label{est:I1FF}
I_1\leq Ch.
\end{align}
Using \fer{def:AF} and \fer{Aq1:split}, we get
\begin{align}
\label{est:ARR4}
& \left|\cA (u_{\tau h}(s,\cdot))(\bx) -\cA_h^\peq{F} (u_{\tau h}(s,\cdot))(\bx)\right| \\
& \hspace{3cm}   \leq 
 \left|\cA^\peq{F} (u_{\tau h}(s,\cdot))(\bx) -\cA_h^\peq{F} (u_{\tau h}(s,\cdot))(\bx)\right| + Cq. \nonumber
\end{align}
Since $u_{\tau h}$ is periodic,  we may use the 
circular convolution theorem to express $ \cA_{h}^{\peq{F},i}$ as a convolution. 
The continuous counterpart, $\cA^{\peq{F},i}$, may be also expressed in convolution form, see \fer{def.conv}. 
Therefore, the estimation of \fer{est:ARR4} is reduced to estimate the difference
\begin{align}
\label{defconv}
 w * H^i(u^n)(\bx) - w_h^\peq{F}\circledast H^i(u^n)(\bx).
\end{align}
In the \emph{interior}, i.e., if $\bx\in\O_\bj\cap\O_\peq{FI}^h$ for some $\bj\in \bJ_\peq{F}$, we have, if $\cQ_Q(u^n[\bj])=q_i$, 
\begin{align*}
 w * H^i(u^n)(\bx) & - w_h^\peq{F}\circledast H^i(u^n)(\bx)   \\
& = \sum_{\bk\in\bJ_\peq{F}}A(u^n[\bk]-q_i) \Big(\int_{\O_\bk} w(\bx-\by)d\by - \int_{\O_{\bj-\bk}} w(\bz)d\bz\Big).
\end{align*}
 Using the change of variable $\bz=\hat \bx_\bj -\by$, with $\hat \bx_\bj = (\bj-{\bf 1}/2)h$, we get 
\begin{align*}
 \Big|\int_{\O_\bk} w(\bx-\by)d\by - \int_{\O_{\bj-\bk}} w(\bz)d\bz \Big|\leq  
 \int_{\O_\bk} \big|w(\bx-\by)-w(\hat \bx_\bj -\by)\big|d\by.  
\end{align*}
If $w\in C^1(B)$, using that $\abs{\bx -\hat\bx_\bj}\leq C h$, we deduce 
\begin{align}
\label{est:wjk}
\sum_{\bk\in\bJ_\peq{F}} \Big|\int_{\O_\bk} w(\bx-\by)d\by - \int_{\O_{\bj-\bk}} w(\bz)d\bz \Big|\leq  
 C h \nor{\nabla w}_{L^1(B)}.  
\end{align}
By density, we obtain the same result for $w\in BV(B)$, but with
$\nor{\nabla w}_{L^1(B)}$ replaced by $\nor{w}_{BV(B)}$. 

To get a bound of \fer{defconv} in the \emph{border}, we use the assumption $(H)_F$. If $w\in L^\infty(\O)$ then we reason like in the deduction of \fer{est:u0FB}. Otherwise, if $w$ is compactly supported in $\O$ then, for $h$ small enough, we have $w_h^\peq{F}(\bx) = w(\bx) =0$ for $\bx\in \O_\peq{FS}^h \cup \O_\peq{FC}^h$. From these observations and  recalling that $u_{\tau h}$ is bounded in $L^\infty(Q_T)$, we deduce  from \fer{est:ARR4} and \fer{est:wjk}  
\begin{align}
\label{est:A2RR2}
I_{4}\leq C(h+q).
\end{align}

We conclude like in the previous sections: using estimates  \fer{est:f1},  \fer{est:f2}, \fer{est:time}, \fer{est:I1FF} and \fer{est:A2RR2} in \fer{conv:fv1}, we obtain, for $\tau,h<1$ ,
\begin{align*}
\int_\O \abs{u(t,\bx)  -  u_{\tau h}(t,\bx)} d\bx \leq &  ~ C(\tau + h +q+r)    \\& + C\int_0^t \int_\O \abs{u(s,\bx)  -  u_{\tau h}(s,\bx)} d\bx ds,
\end{align*}
and Gronwall's lemma implies, 
\begin{align}
 \nor{u-u_{\tau h}}_{L^\infty(0,T;L^1(\O))} \leq C(\tau + h+q). \label{eb.fft}
\end{align}

\begin{remark}
For periodic data, i.e. when the identity \fer{def.convFo} is satisfied, we may obtain stronger estimates of  
\fer{defconv} by using the results of Epstein \cite{Epstein2005}, which establish error bounds in the approximation of $\cF$ by $\cF_d$. In particular, his results are valid for piecewise continuous periodic functions.
\end{remark}

\section{Numerical experiments}\label{sec:numerics}

For the next experiments, we always fix the spatial domain as $\O=(0,1)\times(0,1)$ and the final time as $T=1$. 

In the first set of experiments we focus on the comparison of the execution times among the three methods for different choices of the spatial and range kernels, $J$ and $A$, respectively. When an explicit solution is known, we also compare the relative errors of each method. 

In the second set of experiments we investigate the validity of the error bounds obtained for each method, see \fer{eb.ptw}, \fer{eb.rr} and \fer{eb.fft}, in terms of their  relevant parameters, i.e., the time and spatial mesh steps, $\tau$ and $h$, for all the methods, and the range kernel quantization step, $q$, for the \rr~ and \fft~methods. The spatial kernel quantization step, $r$, only present in the \rr~method, is always taken as the maximum allowed by the precision (double) used in the codes. Thus, the approximation of $J$ is similar for all the methods.  

The experiments have been performed in a standard laptop with an i7 processor. The codes are written in \CC~and run on a single core. The execution time is measured in seconds, using the  \texttt{clock()} library. Although this gives a rough measure, its accuracy is enough for our purposes, given the large differences among the measured execution times of the different methods and data sets. The codes are not specially optimized from the programming point of view, but written literally as described in the previous sections. However, notice that the structure of the \fft~method is specially suitable for parallelization since the slides, $\cA_{h}^{\peq{F},i}(u^n)$, that define the approximation may be computed independently, in separated threads. 

\subsection{Execution time}

\begin{figure}[t]
\centering
\includegraphics[width=4.cm,height=4cm]{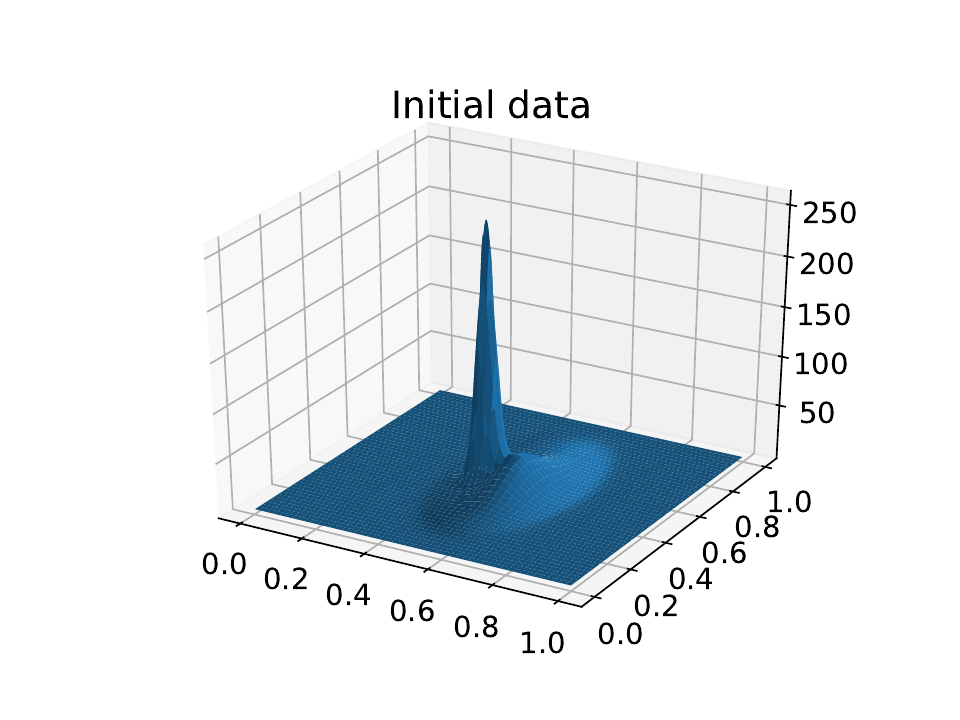}
\includegraphics[width=4.cm,height=4cm]{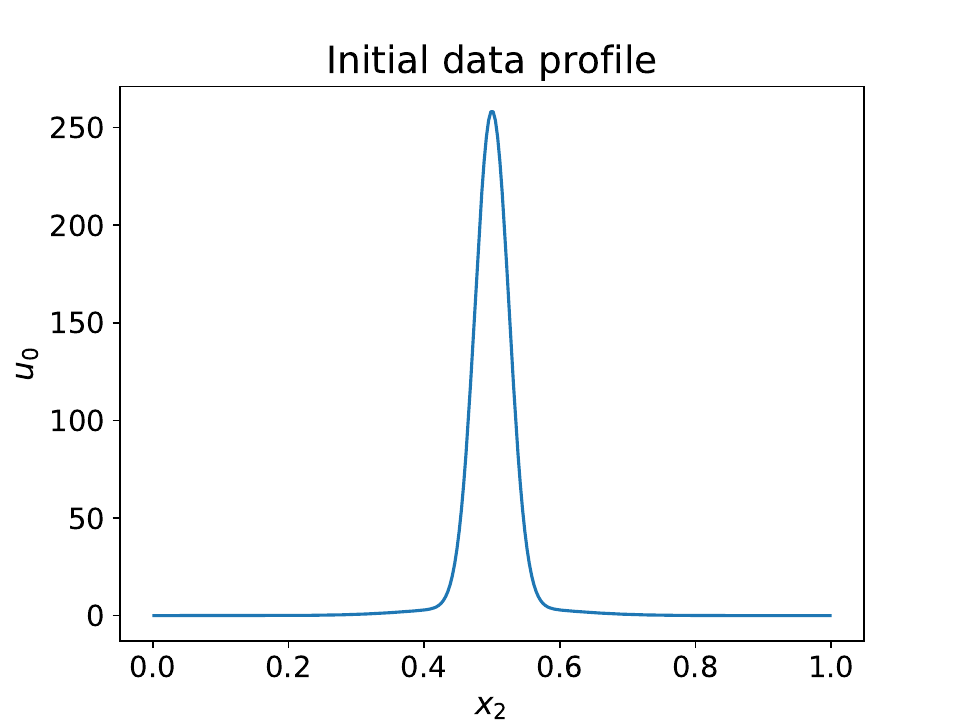}
\includegraphics[width=4.cm,height=4cm]{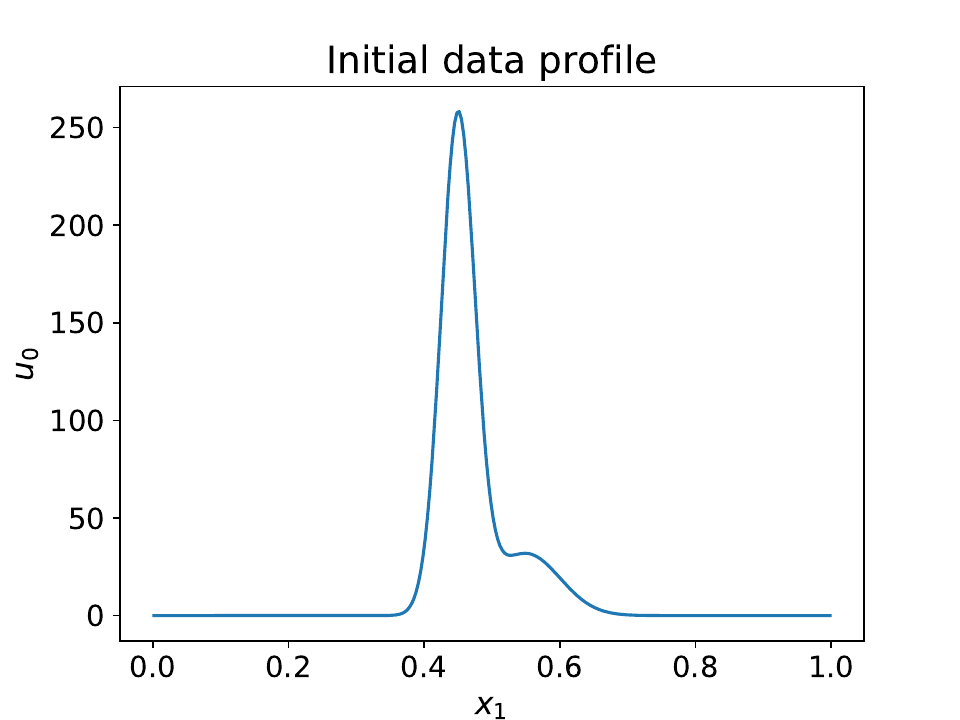}
 \caption{{\small Initial datum and its profiles at $x_1=0.45$ and $x_2=0.5 $.}} 
\label{fig:id}
\end{figure}

The initial datum is given by 
\begin{align}
\label{defu0}
u_0(x_1,x_2) = \sum_{j=1}^2\prod_{i=1}^2 U_0(x_i;\mu_{ji},\sigma_{ji}), 	
\end{align}
where $U_0(x;\mu,\sigma)$ is a Gaussian function of mean $\mu$ and variance $\sigma^2$, with unitary norm in $L^1(\R)$, and 
$\mu_{11} = 0.55, \sigma_{11} = 0.05, \mu_{12} = \mu_{22} = 0.5, \sigma_{12} = 0.1, \mu_{21} = 0.45, \sigma_{21} = \sigma_{22} = 0.025$, see Figure~\ref{fig:id}. 

For the space domain we consider three cases for the uniform mesh \fer{space:disc} or \fer{mesh:fft}, with the number of nodes $J_\O\in\{100,200,300\}$, that is, with step size $h\in\{0.1, 0.005, 0.003 \}$. The time discretization is taken uniform too, with time step $\tau \in \{0.1, 0.01, 0.001\}$. Although in some cases these parameters do not satisfy the sufficient stability condition \fer{cond:est}, 
we always deal with cases in which the schemes are actually convergent.

For the spatial kernel we use either an $L^1$ unitary box spatial kernel of \emph{radius} $r$, 
\begin{align}
\label{def:Jbox}
S(\bx) = \frac{1}{4r^2}1_{B_r(\bx)}	,\quad J(\bx,\by)=S(\bx-\by),
\end{align}
where $B_r(\bx) = (x_1-r,x_1+r)\times(x_2-r,x_2+r)$, with $r\in\{ 0.02,0.05,0.1,0.2,0.3\}$, or a centered and $L^1$ unitary truncated Gaussian kernel, this is,
\begin{align}
\label{def:Jexp}
	G(\bx) = c \exp\Big(-\frac{\abs{\bx}^2}{2\sigma^2}\Big)1_{B_r(\bx)},\quad J(\bx,\by)=G(\bx-\by),
\end{align}
where $c$ is a normalizing constant, and the radius of $B_r(\bx)$ is fixed as $r=6\sigma$, for the data set $\sigma\in\{0.01,0.02,0.03,0.04,0.05\}$.

For the range kernel we take either the identity $A(s)=s$ or $A(s)=\abs{s}^{p-2}s$, for $p=2$ (which is actually the identity) or $p=3$. The latter is the range kernel associated to the nonlocal $p$--Laplacian.

The quantization of the range used in the functional rearrangements and the Fourier based methods is fixed, respectively, as follows. For the former, we take $L_Q=500$ levels of quantization, from $\min_{\bj} u^n[\bj]$ to $\max_{\bj} u^n[\bj]$, where the $n$--dependence indicates that we requantize in each time step. For the latter, we take a smaller number of levels: $10$ in Experiments 1 and 2, and 150 in Experiment 3, without requantization, that is, as a uniform quantization between $\min_{\bj} u^0[\bj]$ and $\max_{\bj} u^0[\bj]$. In the Fourier transform based algorithm, after each time step the approximation is defined as a linear interpolation in terms of the quantized levels: 
\begin{align}
u^{n+1}[\bj] = 	\frac{1}{q}\Big((q_{i+1}-u^{n}[\bj])\cA_{h}^{\peq{F},i} (u^n)[\bj] + (u^{n}[\bj]-q_{i})\cA_{h}^{\peq{F},i+1} (u^n)[\bj]\Big), \label{fft.interp}
\end{align}
see \fer{rangedisc} and \fer{def:cacfq}.

We perform the following experiments:
\begin{enumerate}
\item 	Box spatial kernel and linear range kernel. This is, computationally, the simplest situation. Moreover, an explicit solution may be calculated. Thus, in addition to a time execution comparison of the algorithms we may compute their accuracy. 
\item Box spatial kernel and power-like range kernel. We implement the $p$--Laplacian range kernel $A(s)=\abs{s}^{p-2}s$, but actually taking $p=2$ so that the problem we solve is the same as in Experiment~1. However, it is computationally different since the evaluation of the range kernel $A(u(\by)-u(\bx))$ is performed via an user defined function, instead of directly implemented as $u(\by)-u(\bx)$. This is clearly more time consuming.
\item Gaussian spatial kernel and power-like range kernel. We take the Gaussian space kernel \fer{def:Jbox} and the range kernel $A(s)=\abs{s}^2s$. An explicit solution is not at hand. In addition to the comparison among execution times, we compare the approximations to the finest ($J_\O =300$, $\tau = 0.001$) approximation given by the \ptw~method. In this experiment we introduced a modification of the \fft~algorithm consisting on a zero padding of the samples to the next power of two. As it is well known \cite{fft3}, the computation of the fast Fourier transform is specially efficient in this case. We label this variant as \ffto.

\end{enumerate}

\subsubsection{Experiment 1: Box spatial kernel and linear range kernel}

In this case, we may calculate an explicit solution. Let $r\in(0,1/4)$ and consider the spatial kernel \fer{def:Jexp}.
Let $u_0(\bx)$ be  compactly supported in $(r,1-r)\times(r,1-r)$. Then $u_0(\bx)\int_\O J(\bx-\by)d\by =u_0(\bx)$,
and, for $i,j=1,2$ and $i\neq j$, the convolution term reduces to $I(\bx)= \frac{1}{4r^2}\int_{B_r(\bx)} u_0(\by)d\by$, which may be calculated exactly for our choice of the initial datum, see \fer{defu0}.
Fixing the reaction term as 
\begin{align*}
f(t,\bx,s) = e^{-\lambda t}\big((1-\lambda)u_0(\bx)-I(\bx)\big),
\end{align*}
for some $\lambda\in\R$ ($\lambda = 0.5$ in the experiments), we have that the function $u(t,\bx)=e^{-\lambda t}u_0(\bx)$ is the solution of problem \fer{eq.eq}-\fer{eq.id}. Observe that although the initial datum defined in \fer{defu0} is not compactly supported in $\O$,  its values close to the boundary are so small that its discretized version is actually compactly supported in $\O$.

\begin{figure}[t]
\centering
\includegraphics[width=13cm,height=9cm]{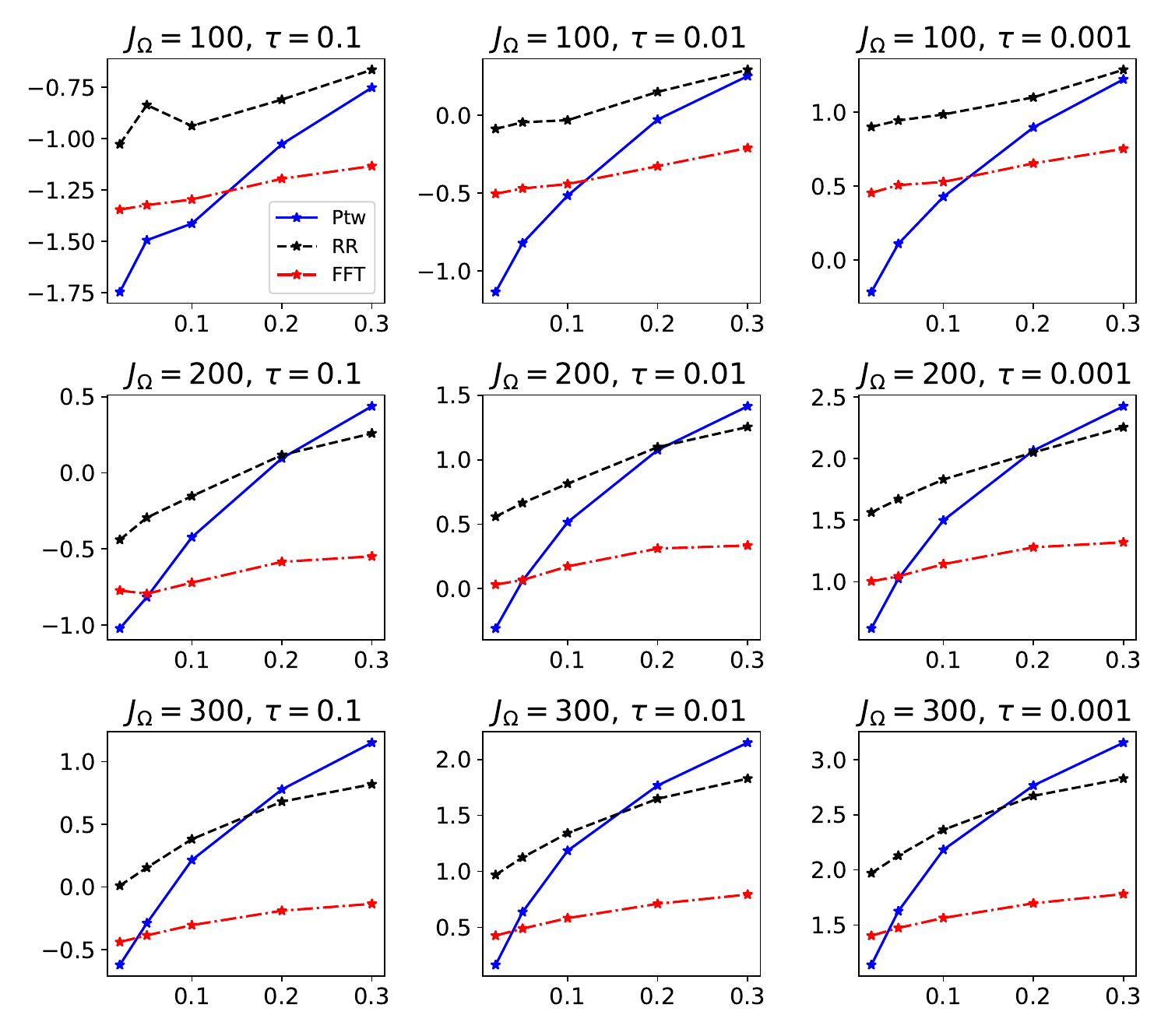}
 \caption{{\small Experiment 1: Execution times (seconds, log$_{10}$ scale) versus radius. Legends are the same for all the plots, see the top-left plot.}} 
\label{fig:exp1et}
\end{figure}

\begin{figure}[t]
\centering
\includegraphics[width=13cm,height=3cm]{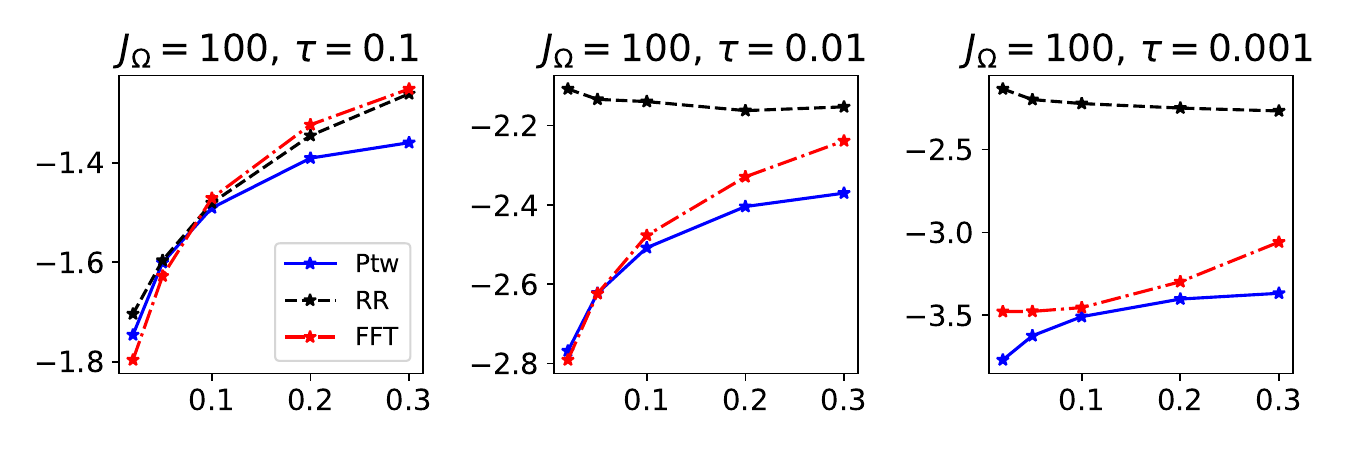}
 \caption{{\small Experiment 1: Relative errors (log$_{10}$ scale) versus radius. We show the case $J_\O = 100$ nodes. The cases $J_\O=200,300$ practically reproduce the same figures.}} 
\label{fig:exp1re}
\end{figure}

In Figures~\ref{fig:exp1et} and \ref{fig:exp1re} we show the execution times and the relative errors, respectively, as functions of the spatial kernel radius, corresponding to each method and the corresponding discretization parameters.
As we shall check in Experiment~2, the linearity of the range kernel is crucial for the \ptw~method to provide lower execution times than the other methods when the radius of the spatial kernel and the number of nodes are small. However, for large radius the \fft~method has always the best figures, and the \rr~method outperforms the \ptw~method. It is also interesting to notice the similarity between the slopes of the execution time versus the radius corresponding to the \fft~and \rr~methods, in contrast to that of the \ptw~method.

With respect to the relative error, computed as 
\begin{align}
RE(u_{\tau h},u) = \sum_{\bj\in J_\O} \abs{u_{\tau h}(T,\bx_\bj)-u(T,\bx_\bj)	} \Big/  \sum_{\bj\in J_\O} \abs{u(T,\bx_\bj)	},\label{def.re}
\end{align}
 the \ptw~discretization is, in general, the most accurate. In any case, the other methods also give good approximations, having the \rr~method a threshold accuracy which depends on the number of quantization levels fixed for the range kernel. This is specially noticeable for small time steps. In the next section, we shall discuss more in detail the error behavior.

\subsubsection{Experiment 2: Box spatial kernel and power-like range kernel}

\begin{figure}[t]
\centering
\includegraphics[width=13cm,height=9cm]{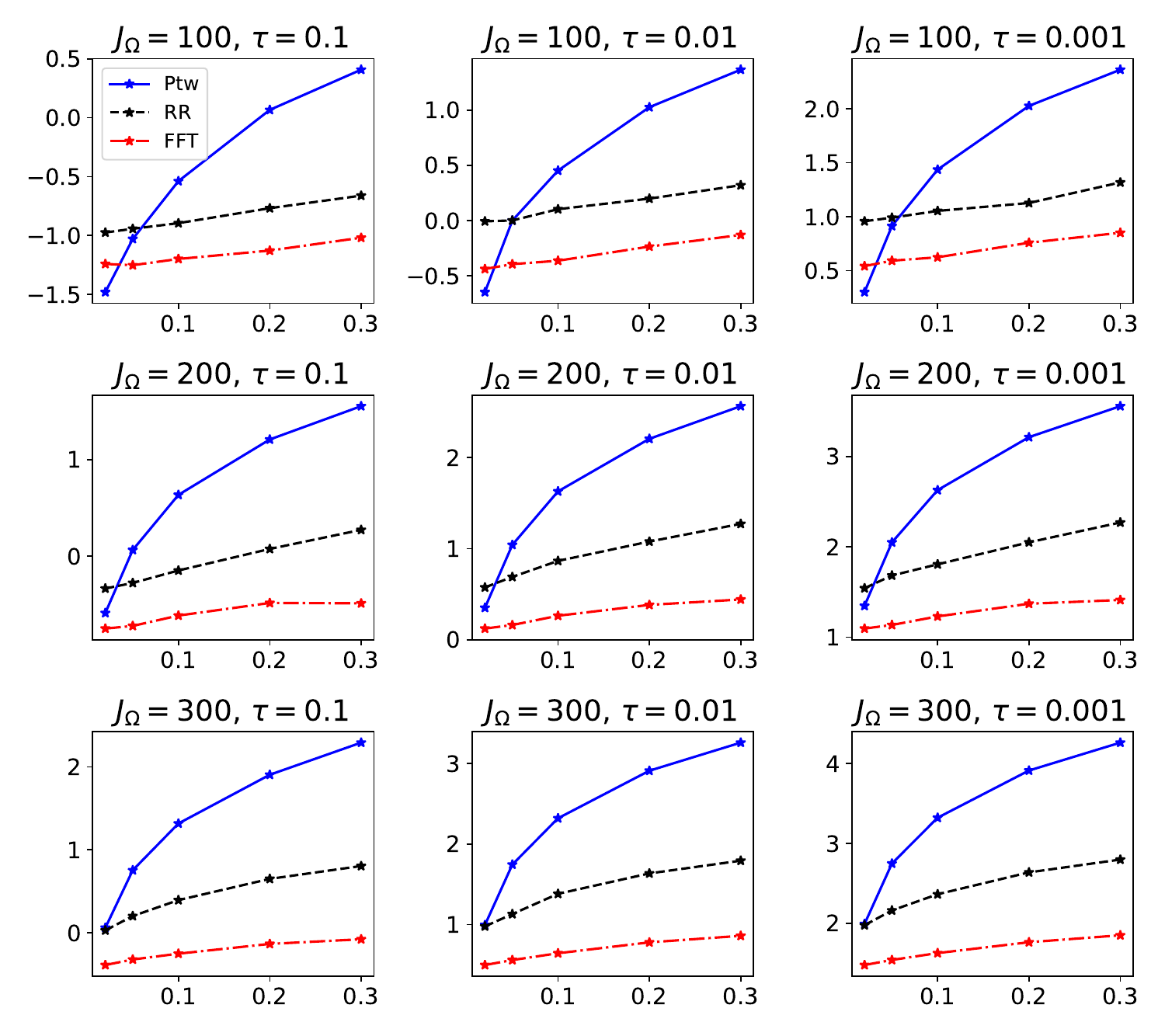}
 \caption{{\small Experiment 2: Execution times (seconds, log$_{10}$ scale) versus radius. Legends are the same for all the plots, see the top-left plot.}} 
\label{fig:exp2et}
\end{figure}

In this experiment we study the increase in execution time due to function evaluation in the range kernel. 
We modify Experiment~1 by replacing the identity range kernel by $A(s)=\abs{s}^{p-2}s$, for $p\geq1$, which defines the $p$--Laplacian nonlocal evolution equation. We, actually, take $p=2$, so that the problem is analytically the same as in Experiment~1. However, the codes are designed to compute the general nonlinear case, which   includes the evaluation of the function $A(s)$, and are therefore slower than the codes of Example~1 that implement this term directly as the identity, without any call to a user defined function.  

In Figure~\ref{fig:exp2et} we show the execution times corresponding to the different methods and discretization parameters. The corresponding relative errors are the same as in Experiment~1. 

The behaviors of the \rr~and the \fft~methods is very similar to those of Experiment~1. However, the \ptw~method suffers from the intensive function evaluation and performs up to  three orders of magnitude slower than its counterparts.  The reason is that for the \rr~method the function evaluations are allocated in a fixed matrix given by $A(q_i-q_j)$  for $i,j=1,\ldots,500$ in our example, see \fer{def:AGmu}, while for the \fft~method only few  \emph{slices} of the range kernel $\cA_{h}^{\peq{F},i} (v_h)[\bj]=\cF_d^{-1}(\cF_d(w_h^\peq{F})\cF_d(H^i(v_h))) [\bj]$, for $i=1,\ldots,10$, are computed, see \fer{def:cacfq}. On the contrary,  the computation of \ptw~in each time step involves two nested loops, one covering all the nodes of the domain, $\bj$, and another covering the support of the spatial kernel (a box of radius $r$), $\bi$, where the evaluation of $A(u^[\bi]-u^[\bj])$ must be accomplished. Therefore, for increasing radius size the computation becomes much more expensive.

\subsubsection{Experiment 3: Gaussian spatial kernel and power-like range kernel}

\begin{figure}[t]
\centering
\includegraphics[width=13cm,height=6cm]{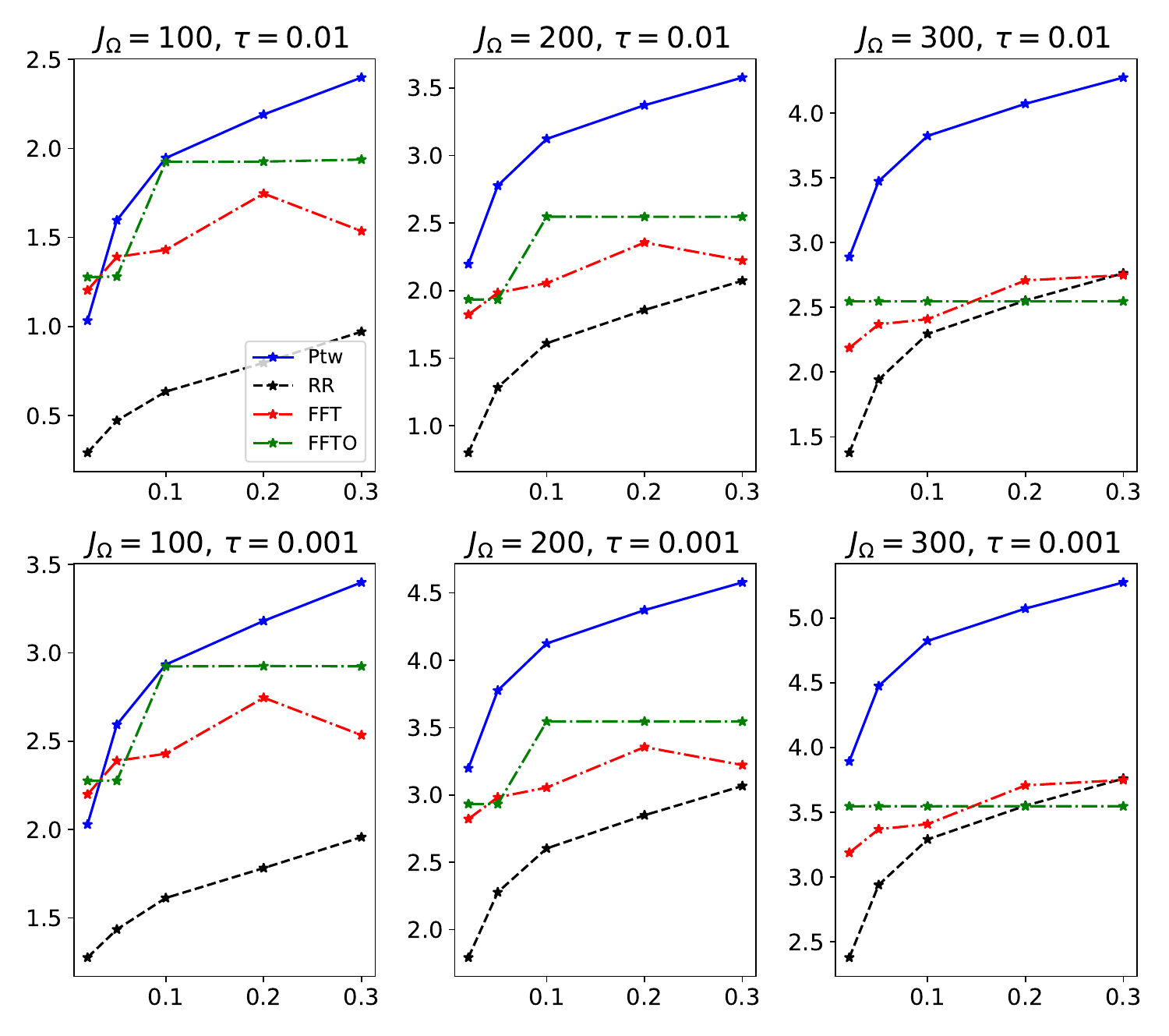}
 \caption{{\small Experiment 3: Execution times (seconds, log$_{10}$ scale) versus radius. Legends are the same for all the plots, see the top-left plot.}} 
\label{fig:exp3et}
\end{figure}

In this experiment we take both the spatial and the range kernels as nonlinear functions defined by the Gaussian \fer{def:Jexp} and the $p$--Laplacian kernel $A(s)=\abs{s}^{p-2}s$, with $p=3$. In addition, we deal with the purely diffusive case, i.e. we set $f\equiv 0$. The stability of the numerical schemes, ensured when condition \fer{cond:est} is satisfied, fails in this case for the time step $\tau=0.1$. Thus, we limit this experiment to the time steps $\tau = 0.01$ and $\tau = 0.001$.  

For some sets of data of the previous experiments, the \fft~method exhibits a somehow erratic behavior in which the execution time decreases or stays nearly flat when the radius of the spatial kernel is increased. See, e.g., the Experiments~1 and 2 with $J_\O=200$. As it is well known, the efficiency of the fast Fourier transform is specially high when the size of the sample to be transformed is a power of two or, more in general, of the form  \cite{fft3} $2^{a_1}3^{a_2}5^{a_3}7^{a_4}11^{a_5}13^{a_6}$, for integer numbers $a_1,\ldots,a_6$. Thus, we incorporated to the set of methods a variant of the \fft~method (\ffto~method) consisting on a slight optimization of the fft by oversampling by zero to the next power of two.  This optimization is effective when the resulting padded samples are not much larger than the original ones.

As we may check in Figure~\ref{fig:exp3et}, only this optimization is capable of outperforming the \rr~method, specially for large spatial kernel radius and for the finer spatial mesh. Let us mention here that the small number (10) of slides (quantized levels) used in the Experiments~1 and 2 for the Fourier transform based method must be replaced by a larger number (150) in order to keep the error in similar bounds than the other methods.

With respect to the \rr~method, the main change in this experiment is the computation of the discretization of the relative rearrangement $J_{*u}$ given by, see \fer{def:AGmu},  
\begin{align}
\label{form.rr}
\sum_{m=0}^{L_\peq{R}}  r_m \mu_m^i[\bj],
\end{align}
where $r_m$ are the discrete level lines of the quantized version of $J$, $\cQ_\peq{R}(J)$, and $\mu_m^i[\bj]$
is the measure of the set of nodes $\{\bk: \cQ_\peq{Q}(u_n[\bk]) = q_i,\text {and } \cQ_\peq{R}(J[\bj,\bk])=r_m \}$. 

In the previous experiments, $J$ has only the two level lines $\{0,1/(4r^2)\}$, making the computation of \fer{form.rr} straightforward. In fact, fast algorithms for this task are easy to implement, see \cite{Porikli2008}. However, if $J$ is a continuous non-constant function, a choice of the number of discrete level lines must be done. Just like for the unknown, $u$, itself.  

In this experiment we fixed, like in the previous, $L_Q=500$ levels of discretization for $u^n$ (i.e., with requantization en each time step). The number $L_R$ of discrete levels of $J$ is fixed as the maximum number of levels given by the precision (double) used in the codes. This is, the number of levels of $G(\bx_k)$ for $\bx_k\in B_r(\bzero)$. This number is  determined by the mesh size and by the radius $r=6\sigma$, where $\sigma\in\{0.01,0.02,0.03,0.04,0.05\}$. For instance, for $J_\O=100$ and $\sigma=0.01$ we get  $L_R=518$, while for $J_\O=300$ and $\sigma=0.05$ we have $L_R=4183$. Of course, both $L_Q$ and $L_R$ are arbitrarily chosen and may be used to find a compromise between accuracy and execution time. 

\begin{figure}[t]
\centering
\includegraphics[width=13cm,height=6cm]{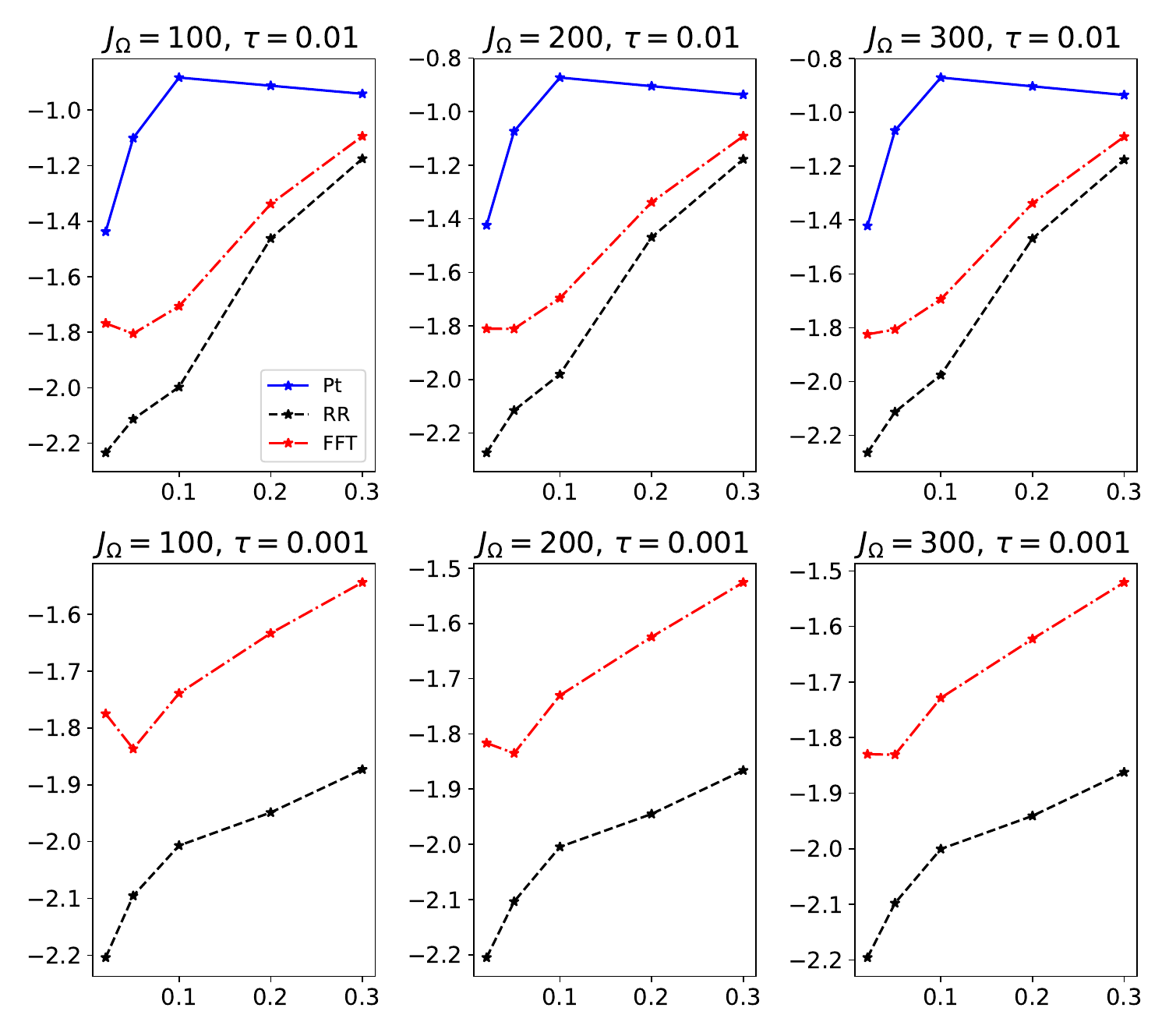}
  \caption{{\small Experiment 3: Relative differences with respect to the \ptw~approximation with $\tau=0.001$ (log$_{10}$ scale) versus radius. Legends are the same for all the plots, see the top-left plot.}} 
\label{fig:exp3re}
\end{figure}

The results of this experiment are certainly clear, see Figures~\ref{fig:exp3et} and \ref{fig:exp3re} . With respect to the execution time, the \ptw~method is much slower than the other methods, specially for large radius sizes, and the \rr~method outperforms the \fft~method up to an order of magnitude for small radius. The only hope for the Fourier transform based method seems to be a further optimization of the sample sizes introduced in the \ffto~implementation.
 
With respect to the error, and taking into account that a exact solution is not known, we have computed the relative differences between the \ptw~approximation obtained for $\tau=0.001$ (assumed to be closer to the exact solution than the others) and the \rr~and the \fft-approximations (or the  \ffto, since they are very close to each other). In contrast with the previous experiments, we see a slight improvement of the \rr~and \fft~approximations for decreasing spatial mesh size. Observe that the only reason for improvement of the approximations because the mesh refinement is due to a higher accuracy in the quadrature of the nonlocal term. However, since this quadrature is reduced to a square of size $4r^2 \ll \abs{\O}$, the improvement is hard to notice.

\subsection{Error}

We consider the range kernel given by \fer{def:Jexp}, for $\sigma\in\{0.01, 0.05\}$, and the linear range kernel $A(s)=s$. Defining the reaction term as 
\begin{align*}
f(t,x) = e^{-\lambda t}\big((m(\bx) - \lambda) u_0(\bx) - J*u_0(\bx)\big),	
\end{align*}
with $ m(\bx) = \int_{\O} J(\bx-\by)d\by$
we have that $u(t,\bx)=e^{-\lambda t}u_0(\bx)$ is an exact solution of problem \fer{eq.eq}-\fer{eq.id}, whatever the initial data $u_0$ we choose. 

The initial datum is taken as $u_0(x_1,x_2)= \Pi_{i=1}^2 U_0(x_i;\mu_0,\sigma_0)$, 	
where $U_0(x;\mu_0,\sigma_0)$ is a Gaussian function of mean $\mu_0=0.5$ and variance $\sigma^2_0=0.01$, and with unitary norm in $L^1(\R)$. Notice that, in practice, since $u_0$ has a small support centered at $(0.5,0.5)$, we have 
$m(\bx)=1$ if $\bx\in\supp(u_0)$.

Due to the linearity of the range kernel and the special choice of $J$ and $u_0$, we may rewrite \fer{eq:evoln} as 
\begin{align*}
u^{n+1}(\bx) = & (1-\lambda \tau) u^n(\bx) +\tau (\lambda -m(\bx))(u^n(\bx)-u(t_n,\bx))\\
&+\tau \big(J*(u^n-u(t_n,\cdot))\big)(\bx).	
\end{align*}
Using any of the error bounds \fer{eb.ptw}, \fer{eb.rr} or \fer{eb.fft}, we deduce an estimate of the type
\begin{align}
 \nor{J*(u^n-u(t_n,\cdot))}_{L^\infty(0,T;L^1(\O))} &\leq \nor{J}_{L^\infty(\O)}	\nor{u^n-u(t_n,\cdot)}_{L^\infty(0,T;L^1(\O))} \nonumber \\
 & \leq C (\tau+h+q+r), \label{b.ju}
\end{align}
implying that the second and third terms of the right hand side of \fer{b.ju} are of second order with respect to the size meshes or, in other words, that 
\begin{align}
u^{n+1}(\bx) = & (1-\lambda \tau) u^n(\bx) +\tau o(\tau+h+q+r),	\label{esexp}
\end{align} 
implementing in this way, up to first order, an explicit Euler approximation of the differential equation  $\p_t u(t,\bx) = \lambda u(t,\bx)$, for each $\bx\in\O$. In particular, the only source of error due to the spatial discretization only contributes to the approximation error as $o(\tau h)$ and pass, in practice, unnoticed in the numerical experiments.

\subsubsection{Experiment 4: Mixed orders of the size meshes}
We run experiments for the following data: time, space and quantization steps $\tau,h,q \in\{0.01,0.005,0.0025\}$ (number of nodes $J_\O= h^{-1}\in\{100,200,400\})$, and for spatial kernel supports determined by $\sigma\in\{0.01,0.05\}$ by the rule $r = 6 \sigma/h  $, where $r$ is the radius of the effective support of the Gaussian. 

The relative error of each experiment, defined in \fer{def.re}, is given in Table~\ref{table:all1}, for the space kernel with smaller support, $\sigma= 0.01$, and in Table~\ref{table:all2}, for the space kernel with larger support, $\sigma= 0.05$. Rows correspond to different space mesh size, $h$, and columns to different time mesh size, $\tau$. For each $(\tau,h)$ we include the errors corresponding to the different quantization size meshes, $q$, implemented in the methods \rr~and \fft. 

For a fair comparison, the \fft~method is implemented without the a posteriori interpolation defined in \fer{fft.interp}. This interpolation is captured by the method labeled as \fft10, for which the corresponding quantization size mesh is multiplied by ten, i.e. $q_{\fftten}\in\{0.1,0.05,0.025\}$.  

Bold numbers in Tables~\ref{table:all1} and ~\ref{table:all2} correspond to the cases in which $\tau=h=q$, for which  the error bounds obtained in \fer{eb.ptw}, \fer{eb.rr} or \fer{eb.fft} predict a relative error of order $o(\tau)$.

Conclusions: As noticed in \fer{esexp}, we may check in both tables the irrelevance of the spatial mesh refinement  for the error figures. We also observe that the relative errors produced by the \rr~and the \fft~methods are always in consonance with the theoretical bounds deduced in   \fer{eb.rr} and \fer{eb.fft}, while the \ptw~method always produces better results, of around one order of magnitude lower than that predicted by \fer{eb.ptw}, most probably due to the smoothness of the solution. This is also observed for the  \rr~ and the \fft~methods when $q$ is reduced, indicating that the main source of error in these methods (and for this example) is the quantization mesh size. Finally, we observe that the relative error corresponding to the \fftten~method is always around one order of magnitude better than the \fft~method, that is, similar to that of the \ptw~method. Indeed, for $\bx\in\O$ we may formally express the nonlocal term as 
\begin{align*}
\cA(v)(\bx)=\int_\R\Phi(s,\bx) \delta(s-v(\bx))ds,
\end{align*}
where $\delta$ is the Dirac delta centered at zero and $\Phi(s,\bx)=\int_\O w(\bx-\by)A(v(\by) - s)d\by$, and approximate it as
\begin{align*}
\cA(v)(\bx)\approx\sum_{i=1}^{L\peq{Q}}\int_{q_i}^{q_{i+1}}\Phi(s,\bx) \rho_n(s-v(\bx))ds,
\end{align*}
 where $\rho_n\to\delta$ as $n\to\infty$. Let $\bx\in\{\by:v(\by)\in  [q_i,q_{i+1})\}$ and $n$ be large enough to have  
 \begin{align*}
 	\cA(v)(\bx) \approx \int_{q_i}^{q_{i+1}}\Phi(s,\bx)ds .
 \end{align*}
 By Taylor's theorem, for $\bx$ fixed we have the first order approximation
\begin{align}
\Phi(s,\bx) &\approx \Phi (q_i,\bx) + \frac{\p\Phi}{\p s}(q_i,\bx)(s-q_i) \nonumber\\
& \approx \Phi (q_i,\bx) + \frac{\Phi(q_{i+1},\bx)-\Phi(q_i,\bx)}{q}(s-q_i). \label{form.in}
\end{align}
 On one hand, using that $(s-q_i)/q = 1 -(q_{i+1}-s)/q$, we deduce that \fer{form.in} is just the interpolation formula defined in \fer{fft.interp} so that the corresponding integral, $\cA (v)(\bx)$ is obtained in the \fftten~method using a trapezoidal (second order) quadrature rule. In the particular case of $A(s)=s$ we have $\p_s\Phi(s,\bx)= const.$, so that the integration of the Taylor's expansion is exact and, therefore, the computation of $\cA(v)(\bx)$ is even of higher order. On the other hand, in the \fft~method we only retain the zero order term of \fer{form.in}, resulting on a first order quadrature rule (left-point), explaining in this way the different orders of approximation of the \fft~and the  \fft10 methods.

\subsubsection{Experiment 5: Same orders of size meshes}
We repite the Experiment~4 for finer meshes of all the discretization parameters, and for a small spatial kernel support ($\sigma=0.01$). Concretely, we take $\tau=h=q\in\{0.005,0.0025,0.0017,0.00125,0.001\}$. Due to the high memory requirements, the \fft~method is unable to run in an standard computer. Thus, we analyze the relative errors of the \ptw, the \rr~and the \fftten~methods. The results are shown in Table~\ref{table:all3}, where the same trend as in Experiment~4 may be observed. In particular, we observe the linear decreasing of the relative error   with respect to the discretization size steps theoretically deduced in \fer{eb.ptw}, \fer{eb.rr} and \fer{eb.fft}.


\begin{table*} 
\small
{\footnotesize
\hspace{-3cm}
\begin{tabular}{|c|c|c|c|}
\hline 
\hline
\multicolumn{4}{|c|}{$\sigma = 0.01$} \\
\hline
$h$
  & $\tau =$ 1.\,e-2 & $\tau =$ 5.\,e-3 & $\tau =$ 2.5\,e-3 \\
\hline
1.\,e-2  & 
\begin{tabular}{|c|c|c|c|}  \multicolumn{4}{c}{\ptw~: \textbf{1.11 \,e-3} } \\ \hline\hline $q$ & \rr~ & \fft~& \fftten~ \\ \hline 1.\,e-2 & \textbf{1.27 \,e-2} & \textbf{1.52 \,e-2}   & \textbf{1.6 \,e-3}   \\ \hline 5.\,e-3 & 7.47 \,e-3    & 8.34 \,e-3   & 1.59 \,e-3   \\ \hline 2.5\,e-3 & 4.22 \,e-3   & 5.02 \,e-3   & 1.62 \,e-3   \\  \end{tabular} &  
\begin{tabular}{|c|c|c|}  \multicolumn{3}{c}{\ptw~: 5.54 \,e-4 } \\ \hline\hline   \rr~ & \fft~& \fftten~\\ \hline  1.29 \,e-2   & 1.44 \,e-2   & 8.7 \,e-4  \\ \hline  7.04 \,e-3   & 7.59 \,e-3   & 7.89 \,e-4  \\ \hline  3.89 \,e-3   & 4.16 \,e-3   & 8.06 \,e-4  \\  \end{tabular} &  
\begin{tabular}{|c|c|c|}  \multicolumn{3}{c}{\ptw~: 2.77 \,e-4 } \\ \hline\hline   \rr~ & \fft~& \fftten~\\ \hline  1.3 \,e-2   & 1.41 \,e-2   & 5.17 \,e-4  \\ \hline  6.91 \,e-3   & 7.13 \,e-3   & 3.98 \,e-4  \\ \hline  3.84 \,e-3   & 3.79 \,e-3   & 4.02 \,e-4  \\  \end{tabular}   \\
\hline
\hline
5.\,e-3  & 
\begin{tabular}{|c|c|c|c|}  \multicolumn{4}{c}{\ptw~: 1.5 \,e-3 } \\ \hline\hline $q$ & \rr~ & \fft~& \fftten~ \\ \hline 1.\,e-2 & 1.11 \,e-2   & 1.51 \,e-2   & 1.6 \,e-3   \\ \hline 5.\,e-3 & 6.4 \,e-3    & 8.3 \,e-3   & 1.6 \,e-3   \\ \hline 2.5\,e-3 & 3.83 \,e-3   & 4.94 \,e-3   & 1.62 \,e-3   \\  \end{tabular} &  
\begin{tabular}{|c|c|c|}  \multicolumn{3}{c}{\ptw~: \textbf{7.49 \,e-4} } \\ \hline\hline   \rr~ & \fft~& \fftten~\\ \hline  1.09 \,e-2   & 1.44 \,e-2   & 8.65 \,e-4  \\ \hline  \textbf{5.91 \,e-3}   & \textbf{7.61 \,e-3}   & \textbf{7.9 \,e-4}  \\ \hline  3.22 \,e-3   & 4.16 \,e-3   & 8.07 \,e-4  \\  \end{tabular} &  
\begin{tabular}{|c|c|c|}  \multicolumn{3}{c}{\ptw~: 3.74 \,e-4 } \\ \hline\hline   \rr~ & \fft~& \fftten~\\ \hline  1.08 \,e-2   & 1.41 \,e-2   & 5.12 \,e-4  \\ \hline  5.75 \,e-3   & 7.16 \,e-3   & 3.97 \,e-4  \\ \hline  3.04 \,e-3   & 3.79 \,e-3   & 4.03 \,e-4  \\  \end{tabular}\\
\hline
\hline
2.5\,e-3  & 
\begin{tabular}{|c|c|c|c|}  \multicolumn{4}{c}{\ptw~: 1.6 \,e-3 } \\ \hline\hline $q$ & \rr~ & \fft~& \fftten~ \\ \hline 1.\,e-2 & 1.06 \,e-2   & 1.51 \,e-2   & 1.6 \,e-3   \\ \hline 5.\,e-3 & 5.97 \,e-3    & 8.28 \,e-3   & 1.6 \,e-3   \\ \hline 2.5\,e-3 & 3.79 \,e-3   & 4.94 \,e-3   & 1.62 \,e-3   \\  \end{tabular} &  
\begin{tabular}{|c|c|c|}  \multicolumn{3}{c}{\ptw~: 7.97 \,e-4 } \\ \hline\hline   \rr~ & \fft~& \fftten~\\ \hline  1.04 \,e-2   & 1.44 \,e-2   & 8.66 \,e-4  \\ \hline  5.51 \,e-3   & 7.56 \,e-3   & 7.91 \,e-4  \\ \hline  3.15 \,e-3   & 4.17 \,e-3   & 8.08 \,e-4  \\  \end{tabular} &  
\begin{tabular}{|c|c|c|}  \multicolumn{3}{c}{\ptw~: \textbf{3.98 \,e-4} } \\ \hline\hline   \rr~ & \fft~& \fftten~\\ \hline  1.03 \,e-2   & 1.41 \,e-2   & 5.1 \,e-4  \\ \hline  5.39 \,e-3   & 7.18 \,e-3   & 3.97 \,e-4  \\ \hline  \textbf{2.91 \,e-3}   & \textbf{3.78 \,e-3}   & \textbf{4.03 \,e-4}  \\  \end{tabular} \\  
\hline
\hline
\end{tabular}
\caption{\small Relative errors corresponding to Experiment 4, $\sigma=0.01$.}
\label{table:all1}
}
\end{table*}



\begin{table*} 
\small
{\footnotesize
\hspace{-3cm}
\begin{tabular}{|c|c|c|c|}
\hline 
\hline
\multicolumn{4}{|c|}{$\sigma_J = 0.05$} \\
\hline
h  & $\tau = $ 1.\,e-2 & $\tau = $ 5.\,e-3 & $\tau = $ 2.5\,e-3 \\
\hline
1.\,e-2  & 
\begin{tabular}{|c|c|c|c|}  \multicolumn{4}{c}{\ptw~: \textbf{1.93 \,e-3} } \\ \hline\hline $q$ & \rr~ & \fft~& \fftten~ \\ \hline 1.\,e-2 & \textbf{8.46 \,e-3}   & \textbf{1.52 \,e-2}   & \textbf{1.57 \,e-3}   \\ \hline 5.\,e-3 & 5.21 \,e-3    & 8.31 \,e-3   & 1.6 \,e-3   \\ \hline 2.5\,e-3 & 2.59 \,e-3   & 5.05 \,e-3   & 1.62 \,e-3   \\  \end{tabular} &  
\begin{tabular}{|c|c|c|}  \multicolumn{3}{c}{\ptw~: 9.62 \,e-4 } \\ \hline\hline   \rr~ & \fft~& \fftten~\\ \hline  8.55 \,e-3   & 1.43 \,e-2   & 7.91 \,e-4  \\ \hline  4.87 \,e-3   & 7.51 \,e-3   & 7.99 \,e-4  \\ \hline  2.34 \,e-3   & 4.15 \,e-3   & 8.08 \,e-4  \\  \end{tabular} &  
\begin{tabular}{|c|c|c|}  \multicolumn{3}{c}{\ptw~: 4.8 \,e-4 } \\ \hline\hline   \rr~ & \fft~& \fftten~\\ \hline  8.62 \,e-3   & 1.4 \,e-2   & 4.33 \,e-4  \\ \hline  4.75 \,e-3   & 7.11 \,e-3   & 3.99 \,e-4  \\ \hline  2.37 \,e-3   & 3.79 \,e-3   & 4.04 \,e-4  \\  \end{tabular} \\
\hline
\hline
5.\,e-3  & 
\begin{tabular}{|c|c|c|c|}  \multicolumn{4}{c}{\ptw~: 1.94 \,e-3 } \\ \hline\hline $q$ & \rr~ & \fft~& \fftten~ \\ \hline 1.\,e-2 & 8.43 \,e-3   & 1.5 \,e-2   & 1.57 \,e-3   \\ \hline 5.\,e-3 & 4.81 \,e-3    & 8.28 \,e-3   & 1.6 \,e-3   \\ \hline 2.5\,e-3 & 2.92 \,e-3   & 4.92 \,e-3   & 1.62 \,e-3   \\  \end{tabular} &  
\begin{tabular}{|c|c|c|}  \multicolumn{3}{c}{\ptw~: \textbf{9.68 \,e-4} } \\ \hline\hline   \rr~ & \fft~& \fftten~\\ \hline  8.38 \,e-3   & 1.43 \,e-2   & 7.89 \,e-4  \\ \hline  \textbf{4.56 \,e-3}   & \textbf{7.53 \,e-3}   & \textbf{8.0 \,e-4}  \\ \hline  2.48 \,e-3   & 4.15 \,e-3   & 8.09 \,e-4  \\  \end{tabular} &  
\begin{tabular}{|c|c|c|}  \multicolumn{3}{c}{\ptw~: 4.84 \,e-4 } \\ \hline\hline   \rr~ & \fft~& \fftten~\\ \hline  8.4 \,e-3   & 1.4 \,e-2   & 4.31 \,e-4  \\ \hline  4.48 \,e-3   & 7.13 \,e-3   & 3.99 \,e-4  \\ \hline  2.39 \,e-3   & 3.79 \,e-3   & 4.05 \,e-4  \\  \end{tabular} \\
\hline
\hline
2.5\,e-3  & 
\begin{tabular}{|c|c|c|c|}  \multicolumn{4}{c}{\ptw~: 1.95 \,e-3 } \\ \hline\hline $q$ & \rr~ & \fft~& \fftten~ \\ \hline 1.\,e-2 & 8.38 \,e-3   & 1.5 \,e-2   & 1.57 \,e-3   \\ \hline 5.\,e-3 & 4.67 \,e-3    & 8.25 \,e-3   & 1.6 \,e-3   \\ \hline 2.5\,e-3 & 3.01 \,e-3   & 4.93 \,e-3   & 1.62 \,e-3   \\  \end{tabular} &  
\begin{tabular}{|c|c|c|}  \multicolumn{3}{c}{\ptw~: 9.7 \,e-4 } \\ \hline\hline   \rr~ & \fft~& \fftten~\\ \hline  8.31 \,e-3   & 1.43 \,e-2   & 7.89 \,e-4  \\ \hline  4.46 \,e-3   & 7.52 \,e-3   & 8.0 \,e-4  \\ \hline  2.52 \,e-3   & 4.16 \,e-3   & 8.1 \,e-4  \\  \end{tabular} &  
\begin{tabular}{|c|c|c|}  \multicolumn{3}{c}{\ptw~: \textbf{4.84 \,e-4} } \\ \hline\hline   \rr~ & \fft~& \fftten~\\ \hline  8.31 \,e-3   & 1.4 \,e-2   & 4.29 \,e-4  \\ \hline  4.43 \,e-3   & 7.17 \,e-3   & 3.99 \,e-4  \\ \hline  \textbf{2.39 \,e-3}   & \textbf{3.77 \,e-3}   & \textbf{4.05 \,e-4}  \\  \end{tabular}  \\  
\hline
\hline
\end{tabular}
\caption{\small Relative errors corresponding to Experiment 4, $\sigma=0.05$}
\label{table:all2}
}
\end{table*}


\newpage

\begin{table*} 
\small
{
\centering
\begin{tabular}{|c|c|c|c|c|c|}
\hline 
\hline
\multicolumn{6}{|c|}{$\sigma_J = 0.01$} \\ 
\hline
 Order expected& 5.\,e-3 & 2.5\,e-3 & 1.7\,e-3 & 1.25 \,e-3& 1. \,e-3\\
\hline
\ptw~ & 7.49 \,e-4  & 3.98 \,e-4  & 4.02 \,e-4  & 2.02 \,e-4  & 1.62 \,e-4 \\ \hline
\rr~ & 5.91 \,e-3  & 2.91 \,e-3  & 2.48 \,e-3  & 1.5 \,e-3  & 1.22 \,e-3   \\ \hline
\fftten~ & 7.9 \,e-4  & 4.03 \,e-4  & 9.4 \,e-4  & 2.03 \,e-4  & 1.62 \,e-4   \\
\hline
\hline
\end{tabular}
\caption{\small Relative errors corresponding to Experiment 5.}
\label{table:all3}
}
\end{table*}

\bigskip 

\no\textbf{References}


\end{document}